\documentclass[12pt,fleqn, a4]{article}

\usepackage[mathscr]{euscript}
\usepackage{amsmath}

\usepackage{amssymb}
\usepackage{latexsym}
\usepackage{graphics}

\usepackage{graphics}
\usepackage{color}
\newcommand{\colval}{0.3}
\definecolor{colone}{gray}{\colval}

\newcommand{\dcb}{\begin{array}{lll}}
\newcommand{\dce}{\end{array}}
\newcommand{\ebe}{\begin{enumerate}\setlength{\baselineskip}{13pt}\setlength{\parskip}{5pt}}
\newcommand{\dbe}{\end{enumerate}}

\newcommand{\ibegin}{\begin{itemize}\setlength{\baselineskip}{19pt}\setlength{\parskip}{7pt}}
\newcommand{\iend}{\end{itemize}}
\newcommand{\ok}{\rule{4pt}{6pt}}
\newcommand{\desb}{\begin{description}}
\newcommand{\dese}{\end{description}}

\newtheorem{Thm}{Theorem}[section]
\newtheorem {Cor}[Thm]{Corollary}
\newtheorem {definition}{Definition}[section]
\newtheorem {pro}{Proposition}[Thm]
\newtheorem {Lemma}[Thm]{Lemma}
\newtheorem {rem}{Remark}[section]
\newtheorem {assumption}[definition]{Assumption}%[section]
\newtheorem {condition}[definition]{Condition}%[section]

\newcommand {\bd}{\begin{definition}}
\newcommand {\ed}{\end{definition}}
\newcommand {\bpro}{\begin{pro}}
\newcommand {\epro}{\end{pro}}
\newcommand {\bl}{\begin{Lemma}}
\newcommand {\el}{\end{Lemma}}
\newcommand {\bcor}{\begin{Cor}}
\newcommand {\ecor}{\end{Cor}}
\newcommand {\brem }{\begin{rem} \rm }
\newcommand {\erem }{\end{rem}}
\newcommand{\bethe}{\begin{Thm}}
\newcommand{\ethe}{\end{Thm}}
\newcommand {\bassumption}{\begin{assumption}}
\newcommand {\eassumption}{\end{assumption}}
\newcommand {\bcondition}{\begin{condition}}
\newcommand {\econdition}{\end{condition}}

\newcommand{\transp}{{^\top\!}}

\def \ind{1\!\!1}
\def\cro#1{\langle #1\rangle}

\begin{document}

\begin{center}
\Large
Drift operator in a viable expansion of information flow\footnote{a working version}
\end{center}

\begin{center}
Shiqi Song

{\footnotesize Laboratoire Analyse et Probabilit\'es\\
Universit\'e d'Evry Val D'Essonne, France\\
shiqi.song@univ-evry.fr}
\end{center}

\

\

\begin{center}
Abstract
\end{center}
{\itshape
A triplet $(\mathbb{P},\mathbb{F},S)$ of a probability measure $\mathbb{P}$, of an information flow $\mathbb{F}=(\mathcal{F}_t)_{t\in\mathbb{R}_+}$, and of an $\mathbb{F}$ adapted asset process $S$, is a financial market model, only if it is viable. In this paper we are concerned with the preservation of the market viability, when the information flow $\mathbb{F}$ is replaced by a bigger one $\mathbb{G}=(\mathcal{G}_t)_{t\geq 0}$ with $\mathcal{G}_t\supset\mathcal{F}_t$. Under the assumption of martingale representation property in $(\mathbb{P},\mathbb{F})$, we prove a necessary and sufficient condition for all viable market in $\mathbb{F}$ to remain viable in $\mathbb{G}$.
}

\vspace{2cm}

\textbf{Key words.} Enlargement of filtrations, hypothesis$(H')$, drift operator, martingale representation property, conditional multiplicity, market viability, structure condition, local martingale deflator, no-arbitrage of first kind.

\textbf{MSC class.} 60G07, 60G44, 60G40.

\hspace{2cm}

\section{Introduction}

A financial market is modeled by a triplet $(\mathbb{P},\mathbb{F},S)$ of a probability measure $\mathbb{P}$, of an information flow $\mathbb{F}=(\mathcal{F}_t)_{t\in\mathbb{R}_+}$, and of an $\mathbb{F}$ adapted asset process $S$. The basic requirement about such a model is its viability. (The notion of viability has been defined in \cite{HK1979} for a general economy. It is now used more specifically to signify that the utility maximization problems have solutions, as in \cite{choulli3, K2010nu, Kar, K2012, K2011se, loew}. The viability is closely linked to the absences of arbitrage opportunity (of some kind) as explained in \cite{loew, KC2010} so that the word sometimes is used to signify no-arbitrage condition.) 
There are situations where one should consider the asset process $S$ in an enlarged information flow $\mathbb{G}=(\mathcal{G}_t)_{t\geq 0}$ with $\mathcal{G}_t\supset\mathcal{F}_t$. The viability of the new market $(\mathbb{P},\mathbb{G},S)$ is not guaranteed. The purpose of this paper is to find such conditions that the viability will be maintained despite the expansion of the information flow.

Concretely, we introduce the notion of the full viability on a time horizon $[0,T]$ for information expansions (cf. subsection \ref{viabilitydef} Assumption \ref{fullviabilitydef}). This means that, for any $\mathbb{F}$ special semimartingale asset process $S$, if $(\mathbb{P},\mathbb{F},S)$ is viable, the expansion market $(\mathbb{P},\mathbb{G},S)$ also is viable on $[0,T]$. Under the assumption of martingale representation property in $(\mathbb{P},\mathbb{F})$, we prove that (cf. Theorem \ref{fullviability}) the full viability on $[0,T]$ is equivalent to the following fact: there exist a (multi-dimensional) $\mathbb{G}$ predictable process $\varphi$ and a (multi-dimensional) $\mathbb{F}$ local martingale $N$, such that (1) for any $\mathbb{F}$ local martingale $X$, the expression $\transp\varphi \centerdot[N,X]^{\mathbb{F}\cdot p}$ is well-defined on $[0,T]$ and $X-\transp\varphi \centerdot[N,X]^{\mathbb{F}\cdot p}$ is a $\mathbb{G}$ local martingale on $[0,T]$; (2) the continuous increasing process $\transp{\varphi}(\centerdot[N^c,\transp N^c]){\varphi}$ is finite on $[0,T]$; (3) the jump increasing process 
$
(
\sum_{0<s\leq t}\left(\frac{\transp {\varphi}_{s}\Delta_{s} N}{1+\transp{\varphi}_{s}\Delta_{s} N}
\right)^2
)^{1/2}
$,
$t\in\mathbb{R}_+$,
is $(\mathbb{P},\mathbb{G})$ locally integrable on $[0,T]$.

It is to note that, if no jumps occurs in $\mathbb{F}$, continuous semimartingale calculus gives a quick solution to the viability problem of the information expansion. The situation becomes radically different when jumps occur, especially because we need to compute and to compare the different projections in $\mathbb{F}$ and in $\mathbb{G}$. (The problem is already difficult, even in the case where the filtration does not change. See \cite{KC2010, K2012}). In this paper we come to a satisfactory result in a general jump situation, thanks to a particular property derived from the martingale representation. In fact, when a process $W$ has the martingale representation property, the jump $\Delta W$ of this process can only take a finite number of \texttt{"}predictable\texttt{"} values. We refer to \cite{song-mrp-drift} for a detailed account, where it is called the finite predictable constraint condition (which has a closed link with the notion of multiplicity introduced in \cite{BEKSY}).  

Usually the martingale representation property is mentioned to characterize a specific process (a Brownian motion, for example). But, in this paper, what is relevant is the stochastic basis $(\mathbb{P},\mathbb{F})$ having a martingale representation property, whatever representation process is. One of the fundamental consequences of the finite predictable constraint condition is the possibility to find a finite family of very simply locally bounded mutually \texttt{"}avoiding\texttt{"} processes which have again the martingale representation property. This possibility reduces considerably the computation complexity and gives much clarity to delicate situations.

The viability property is fundamental for financial market modeling. There exists a huge literature (cf. for example, \cite{choulli3, choulli, Fon, imkeller, ing, kabanov, KC2010, K2012,  Kar, loew, Sch1, Sch2, takaoka, song-takaoka}). Recently, there is a particular attention on the viability problem related to expansions of information flow (cf. \cite{AFK,ACDJ,FJS,song-1405}). It is to notice that, however, the most of the works on expansions of information flow follow two specific ideas : the initial enlargement of filtration or the progressive enlargement of filtration (cf. \cite{CJY, J, JY, Pr, mansuyYor} for definition and properties). In this paper, we take the problem in a very different perspective. We obtain general result, without assumption on the way that $\mathbb{G}$ is constructed from $\mathbb{F}$. It is known (cf. \cite{song, song-l-s}) that the initial or progressive enlargement of filtration are particular situations covered by the so-called local solution method. The methodology of this paper does not take part in this category, adding a new element in the arsenal of filtration analysis.

The concept of information is a fascinating, but also a difficult notion, especially when we want to quantify it. The framework of enlargement of filtrations $\mathbb{F}\subset\mathbb{G}$ offers since long a nice laboratory to test the ideas. In general, no common consensus exists how to quantify the difference between two information flows $\mathbb{F}$ and $\mathbb{G}$. The notion of entropy has been used there (see for example \cite{ankirchnerImkeller,Y}). But a more convincing measurement of information should be the drift operator $\Gamma(X)$, i.e. the operator which gives the drift part of the $\mathbb{F}$ local martingale $X$ in $\mathbb{G}$ (cf. Lemma \ref{linearG}).  This observation is strengthened by the result of the present paper. We have seen that, in the case of our paper, the drift operator takes the form $\Gamma(X)=\transp\varphi \centerdot[N,X]^{\mathbb{F}\cdot p}$ for two factor processes $\varphi$ and $N$ (cf. the drift multiplier assumption in Definition \ref{assump1}), and the full viability of the information expansion is completely determined by the size of the positive quantities $\transp\varphi\centerdot[N^c,\transp N^c]\varphi$ and $\frac{1}{1+\transp\varphi\Delta N}$, which have all the appearance of a measure. See \cite{ankirchner} for complementary discussion. See also \cite{JS} for a use of $\Gamma$ in the study of the martingale representation property in $\mathbb{G}$.

\

\section{Notations and vocabulary}\label{vocabulary}

We employ the vocabulary of stochastic calculus as defined in \cite{HWY, Jacodlivre} with the following specifications.

\

\textbf{Probability space and random variables}

A stochastic basis $(\Omega, \mathcal{A},\mathbb{P},\mathbb{F})$ is a quadruplet, where $(\Omega, \mathcal{A},\mathbb{P})$ is a probability space and $\mathbb{F}$ is a filtration of sub-$\sigma$-algebras of $\mathcal{A}$, satisfying the usual conditions. 

The relationships involving random elements are always in the almost sure sense. For a random variable $X$ and a $\sigma$-algebra $\mathcal{F}$, the expression $X\in\mathcal{F}$ means that $X$ is $\mathcal{F}$-measurable. The notation $\mathbf{L}^p(\mathbb{P},\mathcal{F})$ denotes the space of $p$-times $\mathbb{P}$-integrable $\mathcal{F}$-measurable random variables.

\

\textbf{The processes}

The jump process of a c\` adl\`ag process $X$ is denoted by $\Delta X$, whilst the jump at time $t\geq 0$ is denoted by $\Delta_tX$. By definition, $\Delta_0X=0$ for any c\` adl\`ag process $X$. When we call a process $A$ a process having finite variation, we assume automatically that $A$ is c\`adl\`ag. We denote then by $\mathsf{d}A$ the (signed) random measure that $A$ generates. 

An element $v$ in an Euclidean space $\mathbb{R}^d$ ($d\in\mathbb{N}^*$) is considered as a vertical vector. We denote its transposition by $\transp v$. The components of $v$ will be denoted by $v_h, 1\leq h\leq d$. 

We deal with finite family of real processes $X=(X_h)_{1\leq h\leq d}$. It will be considered as $d$-dimensional vertical vector valued process. The value of a component $X_h$ at time $t\geq 0$ will be denoted by $X_{h,t}$. When $X$ is a semimartingale, we denote by $[X,\transp X]$ the $d\times d$-dimensional matrix valued process whose components are $[X_i,X_j]$ for $1\leq i,j\leq d$.

Two multi-dimensional local martingales $X, X'$ are said mutually pathwisely orthogonal, if $[X_{i},X'_{j}]\equiv 0$ for any component $X_{i}$ of $X$ and $X'_{j}$ of $X'$.

\

\textbf{The projections}

With respect to a filtration $\mathbb{F}$, the notation ${^{\mathbb{F}\cdot p}}\bullet$ denotes the predictable projection, and the notation $\bullet^{\mathbb{F}\cdot p}$ denotes the predictable dual projection. 

\

\textbf{The martingales and the semimartingales}

Fix a probability $\mathbb{P}$ and a filtration $\mathbb{F}$. For any $(\mathbb{P},\mathbb{F})$ special semimartingale $X$, we can decompose $X$ in the form (see \cite[Theorem 7.25]{HWY}) :$$
\dcb
X=X_0+X^m+X^v,\
X^m=X^c+X^{da}+X^{di},
\dce
$$
where $X^m$ is the martingale part of $X$ and $X^v$ is the drift part of $X$, $X^c$ is the continuous martingale part, $X^{da}$ is the part of compensated sum of accessible jumps, $X^{di}$ is the part of compensated sum of totally inaccessible jumps. We recall that this decomposition of $X$ depends on the reference probability and the reference filtration. We recall that every part of the decomposition of $X$, except $X_0$, is assumed null at $t=0$.

\

\textbf{The stochastic integrals}

In this paper we employ the notion of stochastic integral only about the predictable processes. The stochastic integral are defined as 0 at $t=0$. We use a point \texttt{"}$\centerdot$\texttt{"} to indicate the integrator process in a stochastic integral. For example, the stochastic integral of a real predictable process ${H}$ with respect to a real semimartingale $Y$ is denoted by ${H}\centerdot Y$, while the expression $\transp{K}(\centerdot[X,\transp X]){K}$ denotes the process$$
\int_0^t \sum_{i=1}^k\sum_{j=1}^k{K}_{i,s}{K}_{j,s} \mathsf{d}[X_i,X_j]_s,\ t\geq 0,
$$
where ${K}$ is a $k$-dimensional predictable process and $X$ is a $k$-dimensional semimartingale. The expression $\transp{K}(\centerdot[X,\transp X]){K}$ respects the matrix product rule. The value at $t\geq 0$ of a stochastic integral will be denoted, for example, by $\transp{K}(\centerdot[X,\transp X]){K}_t$.

The notion of the stochastic integral with respect to a multi-dimensional local martingale $X$ follows \cite{Jacodlivre}. We say that a (multi-dimensional) $\mathbb{F}$ predictable process is integrable with respect to $X$ under the probability $\mathbb{P}$ in the filtration $\mathbb{F}$, if the non decreasing process $\sqrt{\transp{H}(\centerdot[X,\transp X]){H}}$ is $(\mathbb{P},\mathbb{F})$ locally integrable. For such an integrable process ${H}$, the stochastic integral $\transp{H}\centerdot X$ is well-defined and the bracket process of $\transp{H}\centerdot X$ can be computed using \cite[Remarque(4.36) and Proposition(4.68)]{Jacodlivre}. Note that two different predictable processes may produce the same stochastic integral with respect to $X$. In this case, we say that they are in the same equivalent class (related to $X$).

The notion of multi-dimensional stochastic integral is extended to semimartingales. We refer to \cite{JacShi} for details.

\textbf{Caution.}
Some same definitions will be repeated in different parts of the paper to make the lecture easier.

\

\section{Three fundamental concepts}

Three notions play particular roles in this paper.

\subsection{Enlargements of filtrations and Hypothesis$(H')$ }

Let $\mathbb{F}=(\mathcal{F}_t)_{t\geq 0}$ and $\mathbb{G}=(\mathcal{G}_t)_{t\geq 0}$ be two filtrations on a same probability space such that $\mathcal{F}_t\subset\mathcal{G}_t$. We say then that $\mathbb{G}$ is an expansion (or an enlargement) of the filtration $\mathbb{F}$. Let $T$ be a $\mathbb{G}$ stopping time. We introduce the Hypothesis$(H')$ (cf. \cite{CJY, J, JY, Pr, mansuyYor}):   

\bd
(\textbf{Hypothesis$(H')$} on the time horizon $[0,T]$) We say that {Hypothesis$(H')$} holds for the expansion $\mathbb{F}\subset \mathbb{G}$ on the time horizon $[0,T]$ under the probability $\mathbb{P}$, if all $(\mathbb{P},\mathbb{F})$ local martingale is a $(\mathbb{P},\mathbb{G})$ semimartingale on $[0,T]$.
\ed

Whenever Hypothesis$(H')$ holds, the associated drift operator can be defined (cf. \cite{song-mrp-drift}).

\bl\label{linearG}
Suppose hypothesis$(H')$ on $[0,T]$. Then there exists a linear map $\Gamma$ from the space of all $(\mathbb{P},\mathbb{F})$ local martingales into the space of c\`adl\`ag $\mathbb{G}$-predictable processes on $[0,T]$, with finite variation and null at the origin, such that, for any $(\mathbb{P},\mathbb{F})$ local martingale $X$, $\widetilde{X}:=X-\Gamma(X)$ is a $(\mathbb{P},\mathbb{G})$ local martingale on $[0,T]$. Moreover, if $X$ is a $(\mathbb{P},\mathbb{F})$ local martingale and $H$ is an $\mathbb{F}$ predictable $X$-integrable process, then $H$ is $\Gamma(X)$-integrable and $\Gamma(H\centerdot X)=H\centerdot \Gamma(X)$ on $[0,T]$. The operator $\Gamma$ will be called the drift operator.
\el

\

\subsection{The martingale representation property}\label{MRPsection}

Let us fix a stochastic basis $(\Omega, \mathcal{A},\mathbb{P},\mathbb{F})$. We consider a multi-dimensional stochastic process $W$. We say that $W$ has the martingale representation property in the filtration $\mathbb{F}$ under the probability $\mathbb{P}$, if $W$ is a $(\mathbb{P},\mathbb{F})$ local martingale, and if all $(\mathbb{P},\mathbb{F})$ local martingale is a stochastic integral with respect to $W$. We say that the martingale representation property holds in the filtration $\mathbb{F}$ under the probability $\mathbb{P}$, if there exists a local martingale $W$ which possesses the martingale representation property. In this case we call $W$ the representation process.

\subsubsection{The choice of representation process}\label{www}

Recall the result in \cite{song-mrp-drift}. Suppose the martingale representation property in $(\mathbb{P},\mathbb{F})$. Reconstituting the original representation process if necessary, we can find a particular representation process in the form $W=(W',W'',W''')$ (in juxtaposition of three (multi-dimensional) processes), where $W', W''$ denote respectively the processes defined in \cite[Formulas (4) and (5)]{song-mrp-drift} and $W'''$ denote the process $X^\circ$ in \cite[Section 4.5]{song-mrp-drift}. The processes $W',W'',W'''$ are locally bounded $(\mathbb{P},\mathbb{F})$ local martingales; $W'$ is continuous; $W''$ is purely discontinuous with only accessible jump times; $W'''$ is purely discontinuous with only totally inaccessible jump times; the three components $W',W'',W'''$ are mutually pathwisely orthogonal; the components of $W'$ are mutually pathwisely orthogonal; the components of $W'''$ are mutually pathwisely orthogonal. Let $\mathsf{n}',\mathsf{n}'',\mathsf{n}'''$ denote respectively the dimensions of the three components $W',W'',W'''$. We know that, if $d$ denotes the dimension of the original representation process, $\mathsf{n}'=d$ and $\mathsf{n}''=1+d$. (Notice that some components may be null.) Let $H$ be an $\mathbb{F}$ predictable $W$-integrable process. The vector valued process $H$ is naturally cut into three components $(H',H'',H''')$ corresponding to $(W',W'',W''')$. The pathwise orthogonality implies that $H'_h$ is $W'_h$-integrable for $1\leq h\leq d$, $H''$ is $W''$-integrable, and $H'''_h$ is $W'''_h$-integrable for $1\leq h\leq \mathsf{n}'''$.

The finite predictable constraint condition holds (cf. \cite{song-mrp-drift}). There exists a $\mathsf{n}'''$-dimensional $\mathbb{F}$ predictable process $\alpha'''$ such that$$
\Delta W'''_h = \alpha'''_h \ind_{\{\Delta W'''_h\neq 0\}},\
1\leq h\leq \mathsf{n}'''.
$$
Let $(T_n)_{1\leq n<\mathsf{N}^a}$ ($\mathsf{N}^a$ being a finite or infinite integer) be a sequence of strictly positive $(\mathbb{P},\mathbb{F})$ predictable stopping times such that $[T_n]\cap [T_{n'}]=\emptyset$ for $n\neq n'$ and $\{s\geq 0:\Delta_sW''\neq 0\}\subset\cup_{n\geq 1}[T_n]$. For every $1\leq n<\mathsf{N}^a$, there exists (in a general sense) a partition $(A_{n,0},A_{n,1},\ldots,A_{n,d})$ such that $
\mathcal{F}_{T_n}=\mathcal{F}_{T_n-}\vee\sigma(A_{n,0},A_{n,1},A_{n,2},\ldots,A_{n,d})
$
(a finite multiplicity according to \cite{BEKSY}). Denote $p_{n,k}=\mathbb{P}[A_{n,k}|\mathcal{F}_{T_n-}], 0\leq k\leq d$. We have  
$$
W''_k
=
\sum_{n=1}^{\mathsf{N}^a-}\frac{1}{2^n}(\ind_{A_{n,k}}-p_{n,k})\ind_{[T_n,\infty)}.
$$
(Here $\sum_{n=1}^{\mathsf{N}^a-}$ means $\sum_{1\leq n<\mathsf{N}^a}$.) Let us denote by $\mathsf{a}_{n}$ the vector $(\ind_{A_{n,h}})_{0\leq h\leq d}$ and by $p_n$ the vector $(p_{n,h})_{0\leq h\leq d}$, so that $\Delta_{T_n}W''=\frac{1}{2^n}(\mathsf{a}_{n}-p_n)$ on $\{T_{n}<\infty\}$.

\subsubsection{Coefficient in martingale representation}

If the martingale representation property holds, the $(\mathbb{P},\mathbb{F})$ local martingale $X$ takes all the form $\transp H\centerdot W$ for some $W$-integrable predictable process. We call (any version of) the process $H$ the coefficient of $X$ in its martingale representation with respect to the process $W$. This appellation extends naturally to vector valued local martingales.

When we make computation with the martingale representation, we often need to extract information about a particular stopping time from an entire stochastic integral. The following lemma is proved in \cite[Lemma 3.1]{song-mrp-drift}.

\bl\label{single-jump}
Let $R$ be any $\mathbb{F}$ stopping time. Let $\xi\in\mathbf{L}^1(\mathbb{P},\mathcal{F}_{R})$. Let  $H$ denote any coefficient of the $(\mathbb{P},\mathbb{F})$ martingale  $\xi\ind_{[R,\infty)}-(\xi\ind_{[R,\infty)})^{\mathbb{F}\cdot p}$ in its martingale representation with respect to $W$.
\ebe
\item
If $R$ is predictable, the two predictable processes $H$ and $H\ind_{[R]}$ are in the same equivalent class related to $W$, whose value is determined by the equation on $\{R<\infty\}$
$$
\transp H_{R} \Delta_{R}W=
\xi-\mathbb{E}[\xi|\mathcal{F}_{R-}].
$$
\item
If $R$ is totally inaccessible, the process $H$ satisfies the equations  
$$
\transp H_{R} \Delta_{R}W=\xi\ \mbox{ on $\{R<\infty\}$ and } \
\transp H_{R'} \Delta_{R'}W =0 \mbox{ on $\{R'\neq R, R'<\infty\}$},
$$
for any $\mathbb{F}$ stopping time $R'$.
\dbe
\el

\

\subsection{Raw structure condition}\label{SC-def}

Let a stochastic basis $(\Omega, \mathcal{A},\mathbb{P},\mathbb{F})$ be given.

\bd\label{RSC}
Let $R>0$ be an $\mathbb{F}$ stopping time. We say that a multi-dimensional $(\mathbb{P},\mathbb{F})$ special semimartingale $S$ satisfies the raw structure condition in the filtration $\mathbb{F}$ under the probability $\mathbb{P}$ on the time horizon $[0,R]$, if there exists a real $(\mathbb{P},\mathbb{F})$ local martingale $D$ such that, on the time interval $[0,R]$, $D_0=0, \Delta D<1$, $[S^m_{i}, D]^{\mathbb{F}\cdot p}$ exists, and $S^{v}_{i} = [S^m_{i}, D]^{\mathbb{F}\cdot p}$ for all components $S_i$. We will call $D$ a structure connector. 
\ed

\brem
The concept of Definition \ref{RSC} is motivated by Theorem \ref{deflator-connector} below. Its name \texttt{"}raw structure condition\texttt{"} is inspired from the literature. In \cite{choulli, Sch2}, a condition called \texttt{"}structure condition\texttt{"} has been introduced, which has played important role in the study of incomplete market, especially of the minimal martingale measure and of the Follmer-Schweizer decomposition. The structure condition, defined for $d$-dimensional special semimartingales $X=X_{0}+M+A$, has two aspects. On the one hand, the structure condition requires that drift components $A^{i}, 1\leq i\leq d,$ are absolutely continuous with respect to the oblique brackets $\cro{M^i,M^i}$ of martingale components $M^i, 1\leq i\leq d,$ with density functions $\alpha^i$. On the other hand, it imposes the square integrability condition on $M$ and a specific integrability condition on the density functions $\alpha$. But all these requirements are reinterpretations (under square integrability conditions) of the equation $A = \frac{1}{Z^*_{-}}\centerdot\cro{M,Z^*}$ between the drift part $A$, the martingale part $M$ and a strict martingale densities $Z^*$. This equation is precisely the \texttt{"}raw structure condition\texttt{"}. To discriminate between them, we may roughly qualify the situation by saying that the \texttt{"}structure condition\texttt{"} is an expression of solutions of optimization problems, while the \texttt{"}raw structure condition\texttt{"} is an expression of no-arbitrage problem. (\textbf{N.B.} The notations of this remark will not be in use below.)
\ok
\erem

\bd\label{scdef}
Let $R>0$ be an $\mathbb{F}$ stopping time. We call a strictly positive $\mathbb{F}$ adapted real process $X$ with $X_0=1$, a local martingale deflator on the time horizon $[0,R]$ for a (multi-dimensional) $(\mathbb{P},\mathbb{F})$ special semimartingale $S$, if the processes $X$ and $X S$ are $(\mathbb{P},\mathbb{F})$ local martingales on $[0,R]$. 
\ed

We recall that the existence of local martingale deflators and the raw structure condition are conditions equivalent to the no-arbitrage conditions \texttt{NUPBR} and \texttt{NA1} (cf. \cite{KC2010, takaoka}). We know that, when the no-arbitrage condition \texttt{NUPBR} is satisfied, the market is viable, and vice versa.

\bethe\label{deflator-connector}
Let $R>0$ be an $\mathbb{F}$ stopping time. A (multi-dimensional) special semimartingale $S$ possesses a local martingale deflator $X$ in $(\mathbb{P},\mathbb{F})$ on the time horizon $[0,R]$, if and only if $S$ satisfies the raw structure condition on the time horizon $[0,R]$ with a structure connector $D$. In this case, $X=\mathcal{E}(-D)$ on $[0,R]$.
\ethe

\textbf{Proof.} We know that a strictly positive local martingale is always a Dolean-Dade exponential (cf.\cite[Theorem 9.41]{HWY} or \cite{Jacodlivre,choulli2}).  The lemma is a consequence of the integration by parts formula$$
XS = X_0S_0+ S_-\centerdot X + X_-\centerdot S - X_-\centerdot [S,D].
$$ 
In particular, if $X$ and $XS$ are local martingales on $[0,R]$, $[S,D]$ is locally integrable on $[0,R]$. \ok

\

\section{Main results}\label{DMA}

Together with a given stochastic basis $(\Omega,\mathcal{A}, \mathbb{P},\mathbb{F})$, let $\mathbb{G}$ be an enlargement of $\mathbb{F}$. 

\subsection{Drift multiplier assumption and full viability}\label{viabilitydef}

Our study involves the following notions. In this subsection, $T$ denotes a $\mathbb{G}$ stopping time. 

\bd\label{fullviabilitydef}
(\textbf{Full viability} on $[0,T]$) We say that the expansion $\mathbb{F}\subset \mathbb{G}$ is fully viable on $[0,T]$ under $\mathbb{P}$, if, for any $\mathbb{F}$ asset process $S$ (multi-dimensional special semimartingale with strictly positive components) satisfying the raw structure condition in $(\mathbb{P},\mathbb{F})$, the process $S$ satisfies the raw structure condition in the expanded market environment $(\mathbb{P},\mathbb{G})$ on the time horizon $[0,T]$. 
\ed

\brem\label{locallyboundedM}
As indicated in \cite{song-mrp-drift}, the full viability implies that, for any $(\mathbb{P},\mathbb{F})$ locally bounded local martingale $M$, $M$ satisfies the raw structure condition in $(\mathbb{P},\mathbb{G})$ on $[0,T]$.
\erem

\bd\label{assump1} (\textbf{Drift multiplier assumption}) We say that the drift multiplier assumption holds for the expansion $\mathbb{F}\subset \mathbb{G}$ on $[0,T]$ under $\mathbb{P}$, if 
\
\ebe 
\item
Hypothesis$(H')$ is satisfied for the expansion $\mathbb{F}\subset \mathbb{G}$ on the time horizon $[0,T]$ with a drift operator $\Gamma$;

\item
there exist $N=(N_1,\ldots,N_\mathsf{n})$ an $\mathsf{n}$-dimensional $(\mathbb{P},\mathbb{F})$ local martingale, and ${\varphi}$ an $\mathsf{n}$ dimensional $\mathbb{G}$ predictable process such that, for any $(\mathbb{P},\mathbb{F})$ local martingale $X$, $[N,X]^{\mathbb{F}\cdot p}$ exists, ${\varphi}$ is $[N,X]^{\mathbb{F}\cdot p}$-integrable, and $$
\Gamma(X)=\transp{\varphi}\centerdot [N,X]^{\mathbb{F}\cdot p}
$$
on the time horizon $[0,T]$. 
\dbe
The process $N$ will be called the martingale factor and $\varphi$ will be called the integrated factor of the drift operator $\Gamma$.
\ed

We will need frequently the following consequence of the drift multiplier assumptions \ref{assump1}.

\bl\label{A-G-p}
Suppose the drift multiplier assumptions \ref{assump1}. For any $\mathbb{F}$ adapted c\`adl\`ag process $A$ with $(\mathbb{P},\mathbb{F})$ locally integrable variation, we have $$
A^{\mathbb{G}\cdot p}=A^{\mathbb{F}\cdot p}+\Gamma(A-A^{\mathbb{F}\cdot p})=A^{\mathbb{F}\cdot p}+\transp{\varphi}\centerdot[N,A]^{\mathbb{F}\cdot p}
$$
on $[0,T]$. In particular, for $R$ an $\mathbb{F}$ stopping time, for $\xi\in\mathbf{L}^1(\mathbb{P},\mathcal{F}_{R})$, $$
(\xi\ind_{[R,\infty)})^{\mathbb{G}\cdot p}
=
(\xi\ind_{[R,\infty)})^{\mathbb{F}\cdot p}+\transp{\varphi}\centerdot(\Delta_{R}N\xi\ind_{[R,\infty)})^{\mathbb{F}\cdot p}
$$
on $[0,T]$. If $R$ is $\mathbb{F}$ totally inaccessible, $R$ also is $\mathbb{G}$ totally inaccessible on $[0,T]$.
\el

\textbf{Proof.} We can check the result with \cite[Corollary 5.31]{HWY}. \ \ok

\bcondition\label{1+fin}
For any $\mathbb{F}$ predictable stopping time $R$, for any positive random variable $\xi\in\mathcal{F}_R$, we have $\{\mathbb{E}[\xi|\mathcal{G}_{R-}]>0, R\leq T, R<\infty\}=\{\mathbb{E}[\xi|\mathcal{F}_{R-}]>0, R\leq T, R<\infty\}$. 
\econdition

\brem\label{rq:cR}
Clearly, if the random variable $\xi$ is already in $\mathcal{F}_{R-}$ (or if $\mathcal{F}_{R-}=\mathcal{F}_{R}$), the above set equality holds. Hence, a sufficient condition for Condition \ref{1+fin} to be satisfied is that the filtration $\mathbb{F}$ is quasi-left-continuous (cf. \cite[Definition 3.39]{HWY}).
\erem

\

\subsection{The theorems}

The two notions of full viability and the drift multiplier assumption are closely linked. According to \cite{song-mrp-drift}, under the martingale representation property, the full viability on $[0,T]$ implies the drift multiplier assumption. The aim of this paper is to refine the above result to have an exact relationship between the drift multiplier assumption and the full viability of the expanded information flow. We will prove below the two theorems. Let $T$ be a $\mathbb{G}$ stopping time.

\bethe\label{main}
Suppose that $(\mathbb{P},\mathbb{F})$ satisfies the martingale representation property. Suppose the drift multiplier assumption \ref{assump1} (with the factor processes $N$ and $\varphi$) and Condition \ref{1+fin}. Let $S$ be any $(\mathbb{P},\mathbb{F})$ special semimartingale satisfying the raw structure condition in $(\mathbb{P},\mathbb{F})$ with a structure connector $D$. If the process $
\transp{\varphi}\centerdot[N^c,\transp N^c]{\varphi}
$
is finite on $[0,T]$ and if the process $$
\sqrt{\sum_{0<s\leq t\wedge T}\frac{1}{(1+\transp{\varphi}_{s}\Delta_{s} N)^2}\left(
\Delta_{s}D
+
\transp {\varphi}_{s}\Delta_{s} N  
\right)^2},\ t\in\mathbb{R}_+,
$$
is $(\mathbb{P},\mathbb{G})$ locally integrable, then, $S$ satisfies the raw structure condition on $[0,T]$ in $(\mathbb{P},\mathbb{G})$. 
\ethe

\bethe\label{fullviability}
Suppose that $(\mathbb{P},\mathbb{F})$ satisfies the martingale representation property. Then, $\mathbb{G}$ is fully viable on $[0,T]$, if and only if the drift multiplier assumption \ref{assump1} (with the factor processes $N$ and $\varphi$) and Condition \ref{1+fin} are satisfied such that
\begin{equation}\label{fn-sur-fn}
\dcb
\transp{\varphi}(\centerdot[N^c,\transp N^c]){\varphi}\ \mbox{ is a finite process on $[0,T]$ and }\\
\\
\sqrt{\sum_{0<s\leq t\wedge T}\left(\frac{\transp {\varphi}_{s}\Delta_{s} N}{1+\transp{\varphi}_{s}\Delta_{s} N}
\right)^2},\ t\in\mathbb{R}_+,\
\mbox{ is $(\mathbb{P},\mathbb{G})$ locally integrable.}
\dce
\end{equation}

\ethe

\bcor\label{commondeflator}
Under the conditions of the above theorem, there exists a common $(\mathbb{P},\mathbb{G})$ local martingale deflator for all $(\mathbb{P},\mathbb{F})$ local martingales.
\ecor

\

\section{Raw structure condition decomposed under the martingale representation property}

We now begin the proof of Theorem \ref{main} and Theorem \ref{fullviability}. Recall that, when $(\mathbb{P},\mathbb{F})$ possesses the martingale representation property, we can choose the representation process to ease the computations. We suppose in this section the drift multiplier assumption \ref{assump1} and the following one.

\bassumption\label{assump-mrt}
$(\mathbb{P},\mathbb{F})$ satisfies the martingale representation property, with a representation process $W$ of the form $W=(W',W'',W''')$ satisfying the conditions in subsection \ref{www} with respectively the dimensions $d, 1+d, \mathsf{n}'''$. 
\eassumption

Recall the raw structure condition \ref{RSC}. Let $S$ be a multi-dimensional $\mathbb{F}$ asset process satisfying the raw structure condition in $\mathbb{F}$ with an $\mathbb{F}$ structure connector $D$. Set $M:=S^m$ (in $\mathbb{F}$). Let $T$ be a $\mathbb{G}$ stopping time. Under the drift multiplier assumption \ref{assump1}, the $(\mathbb{P},\mathbb{G})$ canonical decomposition of $S$ on $[0,T]$ is given by $$
S=\widetilde{M}+[D,M]^{\mathbb{F}\cdot p}+\transp{\varphi}\centerdot [N,M]^{\mathbb{F}\cdot p}.
$$
The raw structure condition for $S$ in the expanded market environment $(\mathbb{P},\mathbb{G})$ takes the following form : there exists a $\mathbb{G}$ local martingale $Y$ such that $Y_0=0, \Delta Y<1$, $[Y,\widetilde{M}]^{\mathbb{G}\cdot p}$ exists, and 
\begin{equation}\label{structure-condition}
[Y,\widetilde{M}]^{\mathbb{G}\cdot p}=[D,M]^{\mathbb{F}\cdot p}+\transp{\varphi}\centerdot [N,M]^{\mathbb{F}\cdot p}
\end{equation}
on the time horizon $[0,T]$.

Now, we combine the drift multiplier assumption \ref{assump1} with Assumption \ref{assump-mrt}. We consider the following specific raw structure conditions. (recall $\widetilde{X}=X-\Gamma(X)$.) 
\ebe
\item[. ]\textbf{Continuous raw structure condition related to $D$.}
For $1\leq h\leq d$, there exists a $\mathbb{G}$ predictable $\widetilde{W}'_{h}$-integrable process $K'_h$ such that, on the time horizon $[0,T]$,
\begin{equation}\label{structure-condition-c}
K'_h\centerdot[\widetilde{W}'_h, \widetilde{W}'_h]^{\mathbb{G}\cdot p}
=[D, W'_h]^{\mathbb{F}\cdot p}+\transp{\varphi}\centerdot [N, W'_h]^{\mathbb{F}\cdot p}.
\end{equation}
\vspace{13pt}

\item[. ]\textbf{Accessible raw structure condition related to $D$.}
There exists a $\mathbb{G}$ predictable $\widetilde{W}''$-integrable process $K''$  such that $\transp K''\Delta \widetilde{W}''<1$, and, on the time horizon $[0,T]$,
\begin{equation}\label{structure-condition-da}
\transp K''\centerdot[\widetilde{W}'',\transp \widetilde{W}'']^{\mathbb{G}\cdot p}
=[D,\transp {W}'']^{\mathbb{F}\cdot p}+\transp{\varphi}\centerdot [N,\transp {W}'']^{\mathbb{F}\cdot p}.
\end{equation}
\vspace{13pt}

\item[. ]\textbf{Totally inaccessible raw structure condition related to $D$.}
For $1\leq h\leq \mathsf{n}'''$, there exists a $\mathbb{G}$ predictable $\widetilde{W}'''_h$-integrable process $K'''_h$  such that $K'''_h\Delta \widetilde{W}'''_h<1$, and, on the time horizon $[0,T]$,
\begin{equation}\label{structure-condition-di}
 K'''_h\centerdot[\widetilde{W}'''_h,\widetilde{W}'''_h]^{\mathbb{G}\cdot p}
=[D,{W}'''_h]^{\mathbb{F}\cdot p}+\transp{\varphi}\centerdot [N,{W}'''_h]^{\mathbb{F}\cdot p}.
\end{equation}\vspace{17pt}
\dbe
Note that the above conditions assume in particular that all the stochastic integrals exist. Below, we will call the above conditions the raw structure conditions (\ref{structure-condition-c}), (\ref{structure-condition-da}), (\ref{structure-condition-di}). We will also consider (\ref{structure-condition-c}), (\ref{structure-condition-da}), (\ref{structure-condition-di}) as equations for which, we look for solutions $K'_h,K'',K'''_h$.

\bl\label{piece-ensemble}
Suppose the drift multiplier assumption \ref{assump1} with Assumption \ref{assump-mrt}.
Let $S$ be a multi-dimensional $\mathbb{F}$ asset process satisfying the raw structure condition in $\mathbb{F}$ with an $\mathbb{F}$ structure connector $D$. If the group of the conditions (\ref{structure-condition-c}), (\ref{structure-condition-da}), (\ref{structure-condition-di}) related to $D$ are satisfied, the raw structure condition (\ref{structure-condition}) for $S$ in $\mathbb{G}$ is satisfied.
\el

\textbf{Proof.} 
Write the martingale representation of $M:=S^m$ (in $\mathbb{F}$): $$
\dcb
M&=&\transp H'\centerdot W'+\transp H''\centerdot W''+\transp H'''\centerdot W''',
\dce
$$ 
for some $W$-integrable $\mathbb{F}$ predictable processes $(H', H'', H''')$. Let $K'_h,K'',K'''_h$ be the solutions of respectively (\ref{structure-condition-c}), (\ref{structure-condition-da}), (\ref{structure-condition-di}). Set $K':=(K'_h)_{1\leq h\leq d}$, $K''':=(K'''_h)_{1\leq h\leq \mathsf{n}'''}$ and define $$
Y:=\transp K'\centerdot\widetilde{W}'
+\transp K''\centerdot\widetilde{W}''
+\transp K'''\centerdot\widetilde{W}'''.
$$
Note that, with the drift multiplier assumption, $\Gamma(W'')$ has only $\mathbb{F}$ predictable jumps and $\Gamma(W'),\Gamma(W''')$ are continuous so that $\Delta\widetilde{W}'''=\Delta{W}'''$. With the pathwise orthogonality of the processes ${W}'',{W}'''_h, 1\leq h\leq \mathsf{n}''',$ (cf. subsection \ref{www}), we see that $\Delta Y<1$. With the integrability of $H',H'',H'''$ with respect to separately $W',W'',W'''$ (cf. subsection \ref{www}) and \cite[Lemma 2.2]{JS},$$
\widetilde{M}=\transp H'\centerdot \widetilde{W}'+\transp H''\centerdot \widetilde{W}''+\transp H'''\centerdot \widetilde{W}'''.
$$
Therefore,
$$
\dcb
[Y,\transp \widetilde{M}]
&=&
[\transp K'\centerdot\widetilde{W}'
+\transp K''\centerdot\widetilde{W}''
+\transp K'''\centerdot\widetilde{W}''', \ \ H'\centerdot \widetilde{W}'+H''\centerdot \widetilde{W}''+H'''\centerdot \widetilde{W}''']\\

&=&
\sum_{h=1}^{d} K'_h\centerdot[\widetilde{W}'_h,\widetilde{W}'_h]H'_h
+\transp K''\centerdot[\widetilde{W}'',\transp \widetilde{W}'']H''
+\sum_{h=1}^{\mathsf{n}'''} K'''_h\centerdot[\widetilde{W}'''_h, \widetilde{W}'''_h]H'''_h.
\dce
$$
Let $L>0$ be a constant and define $\mathtt{B}=\{|H|\leq L\}$. We can write
$$
\dcb
&&(\ind_{\mathtt{B}}\centerdot[Y,\transp \widetilde{M}])^{\mathbb{G}\cdot p}\\

&=&
\left(
\sum_{h=1}^{d} K'_h\centerdot[\widetilde{W}'_h,\widetilde{W}'_h]H'_h\ind_{\mathtt{B}}
+\transp K''\centerdot[\widetilde{W}'',\transp \widetilde{W}'']H''\ind_{\mathtt{B}}
+\sum_{h=1}^{\mathsf{n}'''} K'''_h\centerdot[\widetilde{W}'''_h, \widetilde{W}'''_h]H'''_h\ind_{\mathtt{B}}
\right)^{\mathbb{G}\cdot p}
\\

&=&

\sum_{h=1}^{d} K'_h\centerdot[\widetilde{W}'_h,\widetilde{W}'_h]^{\mathbb{G}\cdot p}H'_h\ind_{\mathtt{B}}
+\transp K''\centerdot[\widetilde{W}'',\transp \widetilde{W}'']^{\mathbb{G}\cdot p}H''\ind_{\mathtt{B}}
+\sum_{h=1}^{\mathsf{n}'''} K'''_h\centerdot[\widetilde{W}'''_h, \widetilde{W}'''_h]^{\mathbb{G}\cdot p}H'''_h\ind_{\mathtt{B}}
\\

&=&
+\sum_{h=1}^{d}(\centerdot[D,\transp {W}'_h]^{\mathbb{F}\cdot p}H'_h\ind_{\mathtt{B}}+\transp{\varphi}\centerdot [N,\transp {W}'_h]^{\mathbb{F}\cdot p}H'_h\ind_{\mathtt{B}})\\
&&
+\centerdot[D,\transp {W}'']^{\mathbb{F}\cdot p}H''\ind_{\mathtt{B}}+\transp{\varphi}\centerdot [N,\transp {W}'']^{\mathbb{F}\cdot p}H''\ind_{\mathtt{B}}\\
&&
+
\sum_{h=1}^{\mathsf{n}'''}(\centerdot[D,\transp {W}'''_h]^{\mathbb{F}\cdot p}H'''_h\ind_{\mathtt{B}}+\transp{\varphi}\centerdot [N,\transp {W}'''_h]^{\mathbb{F}\cdot p}H'''_h\ind_{\mathtt{B}})\\

&=&
+\ind_{\mathtt{B}}\centerdot[D,\transp M]^{\mathbb{F}\cdot p}+\ind_{\mathtt{B}}\transp{\varphi}\centerdot [N,\transp M]^{\mathbb{F}\cdot p}.
\dce
$$
This formula for any $L>0$ shows firstly that $[Y,\transp \widetilde{M}]^{\mathbb{G}\cdot p}$ exists and then $$
[Y,\transp \widetilde{M}]^{\mathbb{G}\cdot p}
=
[D,\transp M]^{\mathbb{F}\cdot p}+\transp{\varphi}\centerdot [N,\transp M]^{\mathbb{F}\cdot p}.\ \ok
$$

Below we will solve separately the three raw structure conditions (\ref{structure-condition-c}), (\ref{structure-condition-da}), (\ref{structure-condition-di}) related to $D$. 

\

\section{Solution of the continuous raw structure condition}

The drift multiplier assumption \ref{assump1} and Condition \ref{assump-mrt} are in force. The continuous raw structure condition (\ref{structure-condition-c}) has a quick solution. Let$$
\dcb
D = \transp J'\centerdot W'+\transp J''\centerdot W''+\transp J'''\centerdot W''',\\
N = \transp \zeta'\centerdot W'+\transp \zeta''\centerdot W''+\transp \zeta'''\centerdot W''',
\dce
$$
be the martingale representations of $D,N$ in $\mathbb{F}$.

\bethe\label{answer-c}
Suppose the drift multiplier assumption \ref{assump1} and Condition \ref{assump-mrt}.
The continuous raw structure condition (\ref{structure-condition-c}) related to $D$ is satisfied, if and only if the process $
\transp{\varphi}\centerdot[N^c,\transp N^c]\varphi
$
is a finite process on $[0,T]$. In this case, $K'_h=J'_h+(\transp{\varphi}\zeta')_h$ for $1\leq h\leq d$ are particular solutions. 
\ethe

\textbf{Proof.} 
Recall equation (\ref{structure-condition-c}), for $1\leq h\leq d$,	$$
K'_h\centerdot[\widetilde{W}'_h, \widetilde{W}'_h]^{\mathbb{G}-p}
=[D,W'_h]^{\mathbb{F}-p}+\transp{\varphi}\centerdot [N, W'_h]^{\mathbb{F}-p}.
$$
With the continuity, the equation takes another form
$$
\dcb
K'_h\centerdot[{W}'_h, {W}'_h]
&=&J'_h\centerdot [W'_h, W'_h]+(\transp{\varphi}\zeta')_h\centerdot [W'_h,W'_h].
\dce
$$
Hence, if the continuous raw structure condition related to $D$ has a solution $K'$, necessarily $K'_h=J'_h+(\transp{\varphi}\zeta')_h$ almost surely under the random measure $\mathsf{d}[W'_h,W'_h]$ for $1\leq h\leq d$, and $(\transp{\varphi}\zeta')_h$ is $\widetilde{W}'_h$-integrable ($J'_h$ being by assumption $W'_h$-integrable), i.e., $
\transp{\varphi}\centerdot[N^c,\transp N^c]\varphi
$
is a finite process on $[0,T]$.

Conversely, if the process $
\transp{\varphi}\centerdot[N^c,\transp N^c]\varphi
$
is finite on $[0,T]$, define $K'_h=J'_h+(\transp{\varphi}\zeta')_h, 1\leq h\leq d,$ on $[0,T]$. It forms a solution of the continuous raw structure condition related to $D$.
\ok

\

\section{Solution of the accessible raw structure condition}

The drift multiplier assumption \ref{assump1} and Assumption \ref{assump-mrt} are in force in this section. We consider the accessible raw structure condition (\ref{structure-condition-da}).

\subsection{Equations at the stopping times $T_n$}

Recall (cf. subsection \ref{www}) $(T_n)_{1\leq n<\mathsf{N}^a}$ ($\mathsf{N}^a\leq \infty$) a sequence of strictly positive $(\mathbb{P},\mathbb{F})$ predictable stopping times such that $[T_n]\cap [T_{n'}]=\emptyset$ for $n\neq n'$ and $$
\{s\geq 0:\Delta_sW''\neq 0\}\subset\cup_{n\geq 1}[T_n].
$$  
A $\mathbb{G}$ predictable process $K''$ satisfies the equation (\ref{structure-condition-da}) if and only if 
$$%\label{da-GF}
\dcb
&&\transp K''\centerdot[\widetilde{W}'',\transp \widetilde{W}'']^{\mathbb{G}\cdot p}
=
[D,\transp {W}'']^{\mathbb{F}\cdot p}+\transp{\varphi}\centerdot [N,\transp {W}'']^{\mathbb{F}\cdot p}\\

&=&
\transp J''\centerdot[W'',\transp {W}'']^{\mathbb{F}\cdot p}+\transp{\varphi}\zeta''\centerdot [W'',\transp {W}'']^{\mathbb{F}\cdot p}
=
(\transp J''+\transp{\varphi}\zeta'')\centerdot [W'',\transp {W}'']^{\mathbb{F}\cdot p}
\dce
$$
on $[0,T]$. Computing the jumps at $\mathbb{F}$ stopping times $T_n$, we can also say that $K''$ satisfies the equation (\ref{structure-condition-da}) if and only if, for every $1\leq n<\mathsf{N}^a$, on $\{T_n\leq T, T_n<\infty\}$, $K''_{T_n}$ satisfies the equation 
\begin{equation}\label{at-Tn}
\transp K''_{T_n} \mathbb{E}[\Delta_{T_n}\widetilde{W}'' \transp \Delta_{T_n}\widetilde{W}''|\mathcal{G}_{T_n-}]
=(\transp J''+\transp\overline{\varphi}\zeta'')_{T_n}\mathbb{E}[\Delta_{T_n}W'' \transp \Delta_{T_n}W''|\mathcal{F}_{T_n-}] 
\end{equation}
on $\{T_n\leq T, T_n<\infty\}$. (Recall that $\widetilde{W}''$ has no other jumps than that at the $T_n$'s.)

\

\subsection{Conditional expectation at predictable stopping times $T_n$}

For a fixed $1\leq n<\mathsf{N}^a$, let $(A_{n,0},A_{n,1},\ldots,A_{n,d})$ be the partition which satisfies the relation $
\mathcal{F}_{T_n}=\mathcal{F}_{T_n-}\vee\sigma(A_{n,0},A_{n,1},\ldots,A_{n,d})
$
(cf. subsection \ref{www}). Denote $p_{n,h}=\mathbb{P}[A_{n,h}|\mathcal{F}_{T_n-}]$ and $\overline{p}_{n,h}=\mathbb{P}[A_{n,h}|\mathcal{G}_{T_n-}]$ for $0\leq h\leq d$.  Recall that, in the lemmas below, the drift multiplier assumption \ref{assump1} and Assumption \ref{assump-mrt} are assumed.

\bl\label{F-Tn}
We have
\ebe
\item[. ]
For any finite random variable $\xi\in\mathcal{F}_{T_n}$, the conditional expectation $\mathbb{E}[\xi|\mathcal{F}_{T_n-}]$ is well-defined. Let 
$$
\mathfrak{a}_n(\xi)_h=\ind_{\{p_{n,h}>0\}}\frac{1}{p_{n,h}}\mathbb{E}[\ind_{A_{n,h}}\xi|\mathcal{F}_{T_n-}], \ 0\leq h\leq d.
$$
We have $\xi=\sum_{h=0}^d\mathfrak{a}_n(\xi)_h\ind_{A_{n,h}}$.

\item[. ]
Denote the vector valued random variable $n_{n,h}:=\mathfrak{a}_n(\Delta_{T_n}N)_h, 0\leq h\leq d$. We have $$
(1+\transp{\varphi}_{T_n}n_{n,h})p_{n,h}
=
\mathbb{E}[\ind_{A_{n,h}}|\mathcal{G}_{T_n-}]
=\overline{p}_{n,h},
$$
for $0\leq h\leq d$, on $\{T_n\leq T, T_n<\infty\}$.
  
\dbe
\el

\textbf{Proof.}
The first assertion of the lemma is the consequence of the relation $
\mathcal{F}_{T_n}=\mathcal{F}_{T_n-}\vee\sigma(A_{n,0},A_{n,1},\ldots,A_{n,d}).
$
The second assertion follows from a direct computation of $(\ind_{A_{n,h}}\ind_{[T_n,\infty)})^{\mathbb{G}\cdot p}$ using Lemma \ref{A-G-p}.  \ok

\bl\label{1+fn>0}
For $ 1\leq n<\mathsf{N}^a$ we have $
1+\transp{\varphi}_{T_n}\Delta_{T_n} N> 0.
$
\el

\textbf{Proof.} We compute, for $0\leq h\leq d$,$$
\dcb
0\leq\mathbb{E}[\ind_{\{1+\transp{\varphi}_{T_n}\Delta_{T_n} N \leq 0\}}\ind_{A_{n,h}}|\mathcal{G}_{T_n-}]

&=&\ind_{\{1+\transp{\varphi}_{T_n}n_{n,h} \leq 0\}}(1+\transp{\varphi}_{T_n}n_{n,h})p_h\leq 0.
\dce
$$
It follows that $\ind_{\{1+\transp{\varphi}_{T_n}\Delta_{T_n} N \leq 0\}}\ind_{A_{n,h}}=0$ for $0\leq h\leq d$, i.e., $1+\transp{\varphi}_{T_n}\Delta_{T_n} N>0$. \ok

\bl\label{samekernels}
For $ 1\leq n<\mathsf{N}^a$, on $\{T_n\leq T, T_{n}<\infty\}$, suppose the set equality 
\begin{equation}\label{ptop}
\{0\leq h\leq d: p_{n,h}>0\}=\{0\leq h\leq d: \overline{p}_{n,h}>0\}.
\end{equation}
Then, on $\{T_n\leq T, T_{n}<\infty\}$,  the two matrix $\mathbb{E}[\Delta_{T_n}{W}'' \transp \Delta_{T_n}{W}''|\mathcal{F}_{T_n-}]$ and $\mathbb{E}[\Delta_{T_n}\widetilde{W}'' \transp \Delta_{T_n}\widetilde{W}''|\mathcal{G}_{T_n-}]$ have the same kernel space, which is the space $\mathscr{K}_n$ of $a\in\mathbb{R}\times \mathbb{R}^d$ such that $a$ (as function of its components) is constant on the set $\{0\leq h\leq d: p_{n,h}>0\}$, and the same image space $\mathscr{K}^\perp_n$. There exists a matrix valued $\mathcal{G}_{T_n-}$ measurable random variable $\mathsf{G}_n$ such that $\mathbb{E}[\Delta_{T_n}\widetilde{W}\transp \Delta_{T_n}\widetilde{W}|\mathcal{G}_{T_n-}]\mathsf{G}_n$ is the orthogonal projection $\mathsf{P}_n$ onto $\mathscr{K}^\perp_n$. 
\el

\brem\label{eq:77}
The martingale representation property in $\mathbb{F}$ implies that, for any $\mathbb{F}$ predictable stopping time $R$, $\mathcal{F}_{R-}=\mathcal{F}_{R}$ on $\{R<\infty, R\neq T_{n}, 1\leq n<\mathsf{N}^a\}$. Therefore, because of Lemma \ref{F-Tn}, the validity of the set equalities (\ref{ptop}) in Lemma \ref{samekernels} is equivalent to Condition \ref{1+fin}.
\erem

\textbf{Proof.}
Write $\Delta_{T_n}W''_h=\Delta_{T_n}\widetilde{W}''_h+\Delta_{T_n}\Gamma(W''_h)$ and take the conditioning with respect to $\mathcal{G}_{T_n-}$ on $\{T_n\leq T, T_{n}<\infty\}$. We obtain$$
\Delta_{T_n}\Gamma(W''_h)
=
\mathbb{E}[\Delta_{T_n}W''_h|\mathcal{G}_{T_n-}]
=
\mathbb{E}[\frac{1}{2^n}(\ind_{A_{n,h}}-p_{n,h})|\mathcal{G}_{T_n-}]
=
\frac{1}{2^n}(\overline{p}_{n,h} - p_{n,h}),
$$
so that $\Delta_{T_n}\widetilde{W}''_h=\frac{1}{2^n}(\ind_{A_{n,h}}-\overline{p}_{n,h})$. With this in mind, as well as the set equality $\{0\leq h\leq d: p_{n,h}>0\}=\{0\leq h\leq d: \overline{p}_{n,h}>0\}$, we conclude the first assertion of the lemma. The second assertion can be concluded with \cite[Lemma 5.14]{song-mrp-drift}. \ok

\

\subsection{An intermediate result}

Note that $\widetilde{W}''$ is a $\mathbb{G}$ purely discontinuous local martingale. In fact, $\widetilde{W}''=\ind_{\cup_n[T_n]}\centerdot\widetilde{W}''$.

\bl\label{thmprimeprime}
$K''$ satisfies the accessible raw structure condition (\ref{structure-condition-da}) related to $D$, if and only if Condition \ref{1+fin} holds and, for every $1\leq n<\mathsf{N}^a$, on $\{T_n\leq T, T_n<\infty\}$,  
$$
\transp K''_{T_n} \mathsf{P}_n
=(\transp J''+\transp{\varphi}\zeta'')_{T_n}\mathbb{E}[\Delta_{T_n}W''\transp \Delta_{T_n}W''|\mathcal{F}_{T_n-}]\mathsf{G}_n	
$$
and the process $
\sum_{n=1}^{\mathsf{N}^a-}K''_{T_n}\ind_{[T_n]}
$
is $\widetilde{W}''$ integrable on $[0,T]$, i.e.,
\begin{equation}\label{VKI}
\sqrt{
\sum_{n=1}^{\mathsf{N}^a-}\ind_{\{T_n\leq t\wedge T\}}
\left(
(\transp J''+\transp{\varphi}\zeta'')_{T_n}\mathbb{E}[\Delta_{T_n}W\transp \Delta_{T_n}W|\mathcal{F}_{T_n-}]\mathsf{G}_n
\Delta_{T_n}\widetilde{W}''
\right)^2
}
\end{equation}
is $(\mathbb{P},\mathbb{G})$ locally integrable.
\el

\textbf{Proof.}
\textbf{\texttt{"}If\texttt{"} part.} Suppose Condition \ref{1+fin}.. We note then that  the set equality (\ref{ptop}) in Lemma \ref{samekernels} holds on $\{T_n\leq T,  T_{n}<\infty\}$ for every $1\leq n<\mathsf{N}^a$. Let $K''$ be given by the first formula of the lemma. Lemma \ref{samekernels} implies that $K''$ satisfies formula (\ref{at-Tn}), and hence equation (\ref{structure-condition-da}). Note that, on every $\{T_n\leq T,  T_{n}<\infty\}$, $\Delta_{T_n}\widetilde{W}''_h=\frac{1}{2^n}(\ind_{A_{n,h}}-\overline{p}_{n,h}), 0\leq h\leq d$. This implies that, for any $a\in\mathscr{K}_n$, $\transp a \Delta_{T_n}\widetilde{W}''=0$ (noting that $\ind_{A_{n,h}}=0$ if $\overline{p}_{n,h}=0$), i.e., $\Delta_{T_n}\widetilde{W}''\in\mathscr{K}_n^\perp$, which implies$$
\transp K''_{T_n}\Delta_{T_n}\widetilde{W}''=\transp K''_{T_n}\mathsf{P}_n\Delta_{T_n}\widetilde{W}''.
$$
Together with the first formula, it implies that $K''$ is $\widetilde{W}''$ integrable on $[0,T]$ if and only if expression (\ref{VKI}) is $(\mathbb{P},\mathbb{G})$ locally integrable. 

It remains to prove the inequality $\transp K''_{T_n}\Delta_{T_n}\widetilde{W}''<1$ on $\{T_n\leq T,  T_{n}<\infty\}$. The first formula implies formula (\ref{at-Tn}) which implies
$$
\transp K''_{T_n} \mathbb{E}[\Delta_{T_n}\widetilde{W}''(\ind_{A_{n,h}}-\overline{p}_{n,h})|\mathcal{G}_{T_n-}]
=(\transp J''+\transp{\varphi}\zeta'')_{T_n}\mathbb{E}[\Delta_{T_n}W''(\ind_{A_{n,h}}-{p}_{n,h})|\mathcal{F}_{T_n-}], 
$$
for all $0\leq h\leq d$, on $\{T_n\leq T, T_n<\infty\}$, or equivalently,
$$
\dcb
&&\transp K''_{T_n} \mathbb{E}[\Delta_{T_n}\widetilde{W}''\ind_{A_{n,h}}|\mathcal{G}_{T_n-}]
=
(\transp J''+\transp{\varphi}\zeta'')_{T_n}\mathbb{E}[\Delta_{T_n}W''\ind_{A_{n,h}}|\mathcal{F}_{T_n-}],
\dce
$$
because $$
\dcb
\mathbb{E}[\Delta_{T_n}\widetilde{W}''|\mathcal{G}_{T_n-}]=0,\ \
\mathbb{E}[\Delta_{T_n}W''|\mathcal{F}_{T_n-}]=0.
\dce
$$
We denote by $\mathsf{a}_{n}$ the vector $(\ind_{A_{n,h}})_{0\leq h\leq d}$, by $p_n$ the vector $(p_{n,h})_{0\leq h\leq d}$, by $\overline{p}_n$ the vector $(\overline{p}_{n,h})_{0\leq h\leq d}$, to write $$
\Delta_{T_n}W''=\frac{1}{2^n}(\mathsf{a}_{n}-p_n)\ind_{[T_n,\infty)},\
\Delta_{T_n}\widetilde{W}''=\frac{1}{2^n}(\mathsf{a}_{n}-\overline{p}_n)\ind_{[T_n,\infty)}.
$$
For $0\leq h\leq d$, let $d_{n,h}:=\mathfrak{a}_n(\Delta_{T_n}D)_h$, $n_{n,h}:=\mathfrak{a}_n(\Delta_{T_n}N)_h$ to write
$$
\dcb
&&(\transp J''+\transp{\varphi}\zeta'')_{T_n}\mathbb{E}[\Delta_{T_n}W''\ind_{A_{n,h}}|\mathcal{F}_{T_n-}]\\
&=&
\mathbb{E}[\Delta_{T_n}D\ind_{A_{n,h}}|\mathcal{F}_{T_n-}]
+
\transp{\varphi}_{T_n}\mathbb{E}[\Delta_{T_n}N\ind_{A_{n,h}}|\mathcal{F}_{T_n-}]
=
(d_{n,h}+\transp{\varphi}_{T_n}n_{n,h})
p_{n,h}.
\dce
$$
Let $(\epsilon_0,\ldots,\epsilon_d)$ be the canonical basis in $\mathbb{R}\times\mathbb{R}^d$. By Lemma \ref{F-Tn}, 
$$
\transp K''_{T_n} \mathbb{E}[\Delta_{T_n}\widetilde{W}\ind_{A_{n,h}}|\mathcal{G}_{T_n-}]
=
\frac{1}{2^n}\transp K''_{T_n}(\epsilon_{h}-\overline{p}_n) (1+\transp \varphi_{T_n} n_{n,h}){p}_{n,h}.
$$
Putting them together we obtain a new equality for $K''_{T_n}$:
$$
\frac{1}{2^n}\transp K''_{T_n}(\epsilon_{h}-\overline{p}_n) (1+\transp \varphi_{T_n} n_{n,h}){p}_{n,h}
=
(d_{n,h}+\transp{\varphi}_{T_n}n_{n,h})
p_{n,h},
$$
on $\{T_n\leq T,  T_{n}<\infty\}$, so that, because $\Delta_{T_n}D<1$ and $1+\transp{\varphi}_{T_n}\Delta_{T_n} N> 0$ (cf. Lemma \ref{1+fn>0}),$$
\frac{1}{2^n}\transp K''_{T_n}(\epsilon_{h}-\overline{p}_n)
=
\frac{d_{n,h}+\transp{\varphi}_{T_n}n_{n,h}}{1+\transp{\varphi}_{T_n}n_{n,h}}<1,
$$
on $\{p_{n,h}>0\}\cap A_{n,h}=A_{n,h}$. This proves $\transp K''_{T_n}\Delta_{T_n}\widetilde{W}''<1$ on ${A_{n,h}}$ (for every $0\leq h\leq d$). 

\textbf{\texttt{"}Only if\texttt{"} part.} 
Begin with the identity $$
0=\ind_{\{d_{n,h}-1=0\}}(d_{n,h}-1)\ind_{A_{n,h}}=\ind_{\{d_{n,h}-1=0\}}(\Delta_{T_n}D-1)\ind_{A_{n,h}}.
$$
This means that $\ind_{\{d_{n,h}-1=0\}}\ind_{A_{n,h}}=0$, because $\Delta_{T_n}D-1<0$. Taking conditional expectation with respect to $\mathcal{F}_{T_n-}$ we have $\ind_{\{d_{n,h}-1=0\}}p_{n,h}=0$. 
On the other hand, $\Delta_{T_n}D-1<0$ implies $(d_{n,h}-1)p_{n,h}\leq 0$. 
Combining the two properties, we conclude that, on $\{p_{n,h}>0\}$, $d_{n,h}-1<0$. 

For a $K''$ satisfying the accessible raw structure condition (\ref{structure-condition-da}) related to $D$, it satisfies formula (\ref{at-Tn}). The earlier computations show that formula (\ref{at-Tn}) leads to a formula for $(d_{n,h}-1)
p_{n,h}$: 
$$
\frac{1}{2^n}\transp K''_{T_n}(\epsilon_{h}-\overline{p}_n) (1+\transp \varphi_{T_n} n_{n,h}){p}_{n,h}
-
(1+\transp{\varphi}_{T_n}n_{n,h})
p_{n,h},
=
(d_{n,h}-1)
p_{n,h},
$$
on $\{T_n\leq T,  T_{n}<\infty\}$. Note that $(1+\transp{\varphi}_{T_n}n_{n,h})
p_{n,h}=\overline{p}_{n,h}\geq 0$. We conclude that $$
p_{n,h}>0\ \Rightarrow \ (d_{n,h}-1)
p_{n,h} <0 \ \Rightarrow \ \overline{p}_{n,h}=(1+\transp{\varphi}_{T_n}n_{n,h})p_{n,h}>0.
$$
As a consequence, the set equality (\ref{ptop}) in Lemma \ref{samekernels} holds, which implies, on the one hand, Condition \ref{1+fin}, and on the other hand, the conclusions in Lemma \ref{samekernels}. Now, we can repeat the reasoning in the \texttt{"}If\texttt{"} part to achieve the proof of the lemma.  
\ok

\brem\label{1+fn>0bis}
It is interesting to note that the accessible raw structure condition (\ref{structure-condition-da}) related to $D$ implies that $d_{n,h}-1<0$ and $(1+\transp{\varphi}_{T_n}n_{n,h})>0$, whenever $p_{n,h}>0$.
\erem

\

\subsection{The integrability and conclusion}

The integrability condition (\ref{VKI}) in Lemma \ref{thmprimeprime} looks awful. Using Lemma \ref{F-Tn} we now give a pleasant interpretation of the formula (\ref{VKI})	. Recall that $\mathsf{G}_n$ denotes the $\mathcal{G}_{T_n-}$ measurable random matrix which inverses the matrix $\mathbb{E}[\Delta_{T_n}\widetilde{W}''\transp \Delta_{T_n}\widetilde{W}''|\mathcal{G}_{T_n-}]$ on the space $\mathscr{K}^\perp_n$.

\bl\label{referee}
Under the condition of Lemma \ref{samekernels}, we have
$$
\dcb
&&\sum_{n=1}^{\mathsf{N}^a-}\ind_{\{T_n\leq t\wedge T\}}
\left(
\left(\transp J''_{T_n}
+
\transp{\varphi}_{T_n} \zeta''_{T_n}\right)\mathbb{E}[\Delta_{T_n}W'' \transp\Delta_{T_n}W''|\mathcal{F}_{T_n-}]
\mathsf{G}_n\Delta_{T_n}\widetilde{W}''\right)^2\\

&=&\sum_{n=1}^{\mathsf{N}^a-}\ind_{\{T_n\leq t\wedge T\}}\frac{1}{(1+\transp{\varphi}_{T_n}\Delta_{T_n} N)^2}\left(
\Delta_{T_n}D
+
\transp {\varphi}_{T_n}\Delta_{T_n} N 
\right)^2.
\dce
$$
\el

\textbf{Proof.}
Note that, for $0\leq h\leq d$, $W''_h$ is a bounded process with finite variation. $W''_h$ is always a $\mathbb{G}$ special semimartingale whatever hypothesis$(H')$ is valid or not. We denote always by $\widetilde{W}''$ the $\mathbb{G}$ martingale part of $W''$.

Consider the space $\mathtt{E}=\{0,1,2,\ldots,d\}$. We fix $1\leq n<\mathsf{N}^a$ and endow $\mathtt{E}$ with two (random) probability measures $$
\dcb
\mathsf{m}[\{h\}]:=p_{n,h},\\
\overline{\mathsf{m}}[\{h\}]:=\overline{p}_{n,h} = (1+\transp{\varphi}_{T_n}n_{n,h})p_{n,h},\ 0\leq h\leq d.
\dce
$$
Let $\boldsymbol{\epsilon}=(\epsilon_h)_{0\leq h\leq d}$ denote the canonical basis in $\mathbb{R}\times\mathbb{R}^{d}$. Let $d_{n,h}:=\mathfrak{a}_n(\Delta_{T_n}D)_h$, $n_{n,h}:=\mathfrak{a}_n(\Delta_{T_n}N)_h$, and $$
\dcb
&&w_{n,h}:=\mathfrak{a}_n(\Delta_{T_n}W'')_h
=
\frac{1}{2^n}\ind_{\{p_{n,h}>0\}}\left({\epsilon}_h - p_n\right).
\dce
$$ 
(Recall the notations $p_{n}, \overline{p}_{n}, \mathsf{a}_{n}$ in the proof of Lemma \ref{thmprimeprime}.) We define then the function $\mathsf{d}:=\sum_{h=0}^{d}d_{n,h}\ind_{\{h\}}$, $\mathsf{n}:=\sum_{h=0}^{d}n_{n,h}\ind_{\{h\}}$, and $$
\mathsf{w}:
=
\frac{1}{2^n}\sum_{h=0}^{d}{\epsilon}_h\ind_{\{h\}}
-
\frac{1}{2^n}\sum_{h=0}^{d}\ind_{\{h\}}p_n
=
\frac{1}{2^n}\sum_{h=0}^{d}{\epsilon}_h\ind_{\{h\}}
-
\frac{1}{2^n}p_n,
$$ 
on $\mathtt{E}$. As $\mathbb{E}_{{\mathsf{m}}}[\ind_{\{p_{n,h}=0\}}\ind_{\{h\}}]=0$, $\mathsf{w}$ is $\mathsf{m}-a.s.$ and $\overline{\mathsf{m}}-a.s.$ equal to
$$
\sum_{h=0}^{d}w_{n,h}\ind_{\{h\}}
=
\frac{1}{2^n}\sum_{h=0}^{d}\ind_{\{p_{n,h}>0\}}{\epsilon}_h\ind_{\{h\}}
-
\frac{1}{2^n}\sum_{h=0}^{d}\ind_{\{p_{n,h}>0\}}\ind_{\{h\}}p_n.
$$
The function $\mathsf{w}$ is $(1+d)$-dimensional vector valued. We denote by $\mathsf{w}_k$ its $k$th component, which is the real function$$
\mathsf{w}_k
=
\frac{1}{2^n}\ind_{\{k\}}
-
\frac{1}{2^n}p_{n,k},
$$ 
on $\mathtt{E}$. Be careful: do not confound it with $\mathsf{w}(h)$ which is a vector. We have $$
\mathbb{E}_{{\mathsf{m}}}[\mathsf{w}_k]=0,\
\mathbb{E}_{\overline{\mathsf{m}}}[\mathsf{w}_k]=\frac{1}{2^n}(\overline{p}_{n,k}-p_{n,k}),
$$
so that, on $\{T_{n}\leq T, T_n<\infty\}$,$$
\Delta_{T_n}\widetilde{W}''=\frac{1}{2^n}(\mathsf{a}_n-\overline{p}_n)
= \frac{1}{2^n}(\mathsf{a}_n-p_n - (\overline{p}_n - p_n)) 
=
\Delta_{T_n}{W}'' - \mathbb{E}_{\overline{\mathsf{m}}}[\mathsf{w}].
$$
For a function $F$ we compute.$$
\dcb
&&\mathbb{E}[F(\Delta_{T_n}{W}'')|\mathcal{G}_{T_n-}]
=
\sum_{h=0}^{d}F(w_{n,h})\mathbb{E}[\ind_{A_{n,h}}|\mathcal{G}_{T_n-}]

=
\mathbb{E}_{\overline{\mathsf{m}}}[F(\mathsf{w})].
\dce
$$
Similarly, 
$\mathbb{E}[F(\Delta_{T_n}{W}'')|\mathcal{F}_{T_n-}]
=\mathbb{E}_{{\mathsf{m}}}[F(\mathsf{w})].
$
Let $$
\dcb
x:&=&\left(J''_{T_n}
+
\transp\zeta''_{T_n}{\varphi}_{T_n} \right),\\
y:&=&\mathsf{G}_n \mathbb{E}[\Delta_{T_n}W'' \transp\Delta_{T_n}W''|\mathcal{F}_{T_n-}] x.
\dce
$$
Then, for all $z\in\mathbb{R}\times\mathbb{R}^{d}$, we can write
\begin{equation}\label{zWWyx}
\transp z \mathbb{E}[\Delta_{T_n}\widetilde{W}''\transp \Delta_{T_n}\widetilde{W}''|\mathcal{G}_{T_n-}] y
=
\transp z \mathbb{E}[\Delta_{T_n}{W}''\transp \Delta_{T_n}{W}''|\mathcal{F}_{T_n-}] x,
\end{equation} 
because $\mathbb{E}[\Delta_{T_n}{W}''\transp \Delta_{T_n}{W}''|\mathcal{F}_{T_n-}] x\in \mathscr{K}^\perp_n$. As$$
\mathbb{E}[\Delta_{T_n}\widetilde{W}''\transp \Delta_{T_n}\widetilde{W}''|\mathcal{G}_{T_n-}]
=
\mathbb{E}[\Delta_{T_n}{W}''\transp \Delta_{T_n}\widetilde{W}''|\mathcal{G}_{T_n-}],
$$
the equation (\ref{zWWyx}) becomes
$$
\dcb
\mathbb{E}[(\transp z\Delta_{T_n}{W}'')(\transp \Delta_{T_n}\widetilde{W}'' y)|\mathcal{G}_{T_n-}]
=
\mathbb{E}[(\transp z \Delta_{T_n}{W}'')(\transp \Delta_{T_n}{W}''x)|\mathcal{F}_{T_n-}],
\dce
$$
or equivalently
$$
\dcb
\mathbb{E}_{\overline{\mathsf{m}}}[(\transp z\mathsf{w})\transp(\mathsf{w}-\mathbb{E}_{\overline{\mathsf{m}}}[\mathsf{w}]) y]
=
\mathbb{E}_{{\mathsf{m}}}[(\transp z \mathsf{w})(\transp \mathsf{w}x)].
\dce
$$
Set the function $\mathsf{q}:=\sum_{h=0}^{\mathsf{n}''}(1+\transp{\varphi}_{T_n}n_{n,h})\ind_{\{h\}}$ on $\mathtt{E}$. We note that $\mathsf{q}=\frac{\mathsf{d}\overline{\mathsf{m}}}{\mathsf{d}{\mathsf{m}}}$. We have therefore
$$
\dcb
\mathbb{E}_{{\mathsf{m}}}[(\transp z\mathsf{w})\mathsf{q}(\transp\mathsf{w}-\mathbb{E}_{\overline{\mathsf{m}}}[\transp\mathsf{w}]) y]
=
\mathbb{E}_{{\mathsf{m}}}[(\transp z \mathsf{w})(\transp \mathsf{w}x)].
\dce
$$
For any vector $a=(a_0,a_1,\ldots,a_d)\in\mathbb{R}\times \mathbb{R}^d$ such that $\transp a p_n=0$, we have$$
\transp a \mathsf{w}
=
\sum_{k=0}^da_k\mathsf{w}_k
=
\sum_{k=0}^da_k(\frac{1}{2^n}\ind_{\{k\}}
-
\frac{1}{2^n}p_{n,k})
=
\frac{1}{2^n}\sum_{k=0}^da_k\ind_{\{k\}}.
$$
This means that the functions of the form $(\transp a \mathsf{w})$ generate the space of all functions on $\mathtt{E}$ with null ${\mathsf{m}}$-expectation. But,$$
\mathbb{E}_{{\mathsf{m}}}[\mathsf{q}(\transp\mathsf{w}-\mathbb{E}{\overline{\mathsf{m}}}[\transp\mathsf{w}]) y]
=
\mathbb{E}_{\overline{\mathsf{m}}}[(\transp\mathsf{w}-\mathbb{E}_{\overline{\mathsf{m}}}[\transp\mathsf{w}]) y]
=0.
$$
Hence, the above identity for all $\transp z\mathsf{w}$ implies$$
\mathsf{q}(\transp\mathsf{w}-\mathbb{E}_{\overline{\mathsf{m}}}[\transp\mathsf{w}]) y
=
\transp \mathsf{w}x\  \mbox{ or  }\
(\transp\mathsf{w}-\mathbb{E}_{\overline{\mathsf{m}}}[\transp\mathsf{w}]) y
=
\frac{1}{\mathsf{q}}\transp \mathsf{w}x,
\
\mathsf{m}-a.s.,
$$
(cf. Remark \ref{1+fn>0bis}). Regarding the values at every $0\leq h\leq d$ with $p_{n,h}>0$,
$$
(\transp w_{n,h}-\mathbb{E}_{\overline{\mathsf{m}}}[\transp\mathsf{w}]) y
=
\frac{1}{(1+\transp{\varphi}_{T_n}n_k)}\transp w_{n,h} x.
$$
Consider now the process
$$
\dcb
\sum_{n=1}^{\mathsf{N}^a-}\ind_{\{T_n\leq t\wedge T\}}
\left(
\left(\transp J''_{T_n}
+
\transp{\varphi}_{T_n} \zeta''_{T_n}\right)\mathbb{E}[\Delta_{T_n}W'' \transp\Delta_{T_n}W''|\mathcal{F}_{T_n-}]
\mathsf{G}_n\Delta_{T_n}\widetilde{W}''\right)^2.
\dce
$$
We have$$
\dcb
&&
\left(\transp J''_{T_n}
+
\transp{\varphi}_{T_n} \zeta''_{T_n}\right)\mathbb{E}[\Delta_{T_n}W'' \transp\Delta_{T_n}W''|\mathcal{F}_{T_n-}]
\mathsf{G}_n\Delta_{T_n}\widetilde{W}''\\

&=&
\transp y \Delta_{T_n}\widetilde{W}''
=
\sum_{h=0}^{d}\transp\Delta_{T_n}\widetilde{W}'' y\ind_{A_{n,h}}
=
\sum_{h=0}^{d}\transp (w_{n,h}-\mathbb{E}_{\overline{\mathsf{m}}}[\mathsf{w}])y\ind_{A_{n,h}}\\

&=&
\sum_{h=0}^{d}\frac{1}{(1+\transp{\varphi}_{T_n}n_k)}\transp w_{n,h} x\ind_{A_{n,h}}.
\dce
$$
It implies
$$
\dcb
&&
\left(
\left(\transp J''_{T_n}
+
\transp{\varphi}_{T_n} \zeta''_{T_n}\right)\mathbb{E}[\Delta_{T_n}W'' \transp\Delta_{T_n}W''|\mathcal{F}_{T_n-}]
\mathsf{G}_n\Delta_{T_n}\widetilde{W}''
\right)^2\\

&=&
\left(
\sum_{h=0}^{d}\frac{1}{(1+\transp{\varphi}_{T_n}n_k)}\transp w_{n,h} x\ind_{A_{n,h}}
\right)^2

=
\sum_{h=0}^{d}\frac{1}{(1+\transp{\varphi}_{T_n}n_k)^2}\left(
\transp w_{n,h} x 
\right)^2\ind_{A_{n,h}}\\

&=&
\sum_{h=0}^{d}\frac{1}{(1+\transp{\varphi}_{T_n}\Delta_{T_n} N)^2}\left(
\transp x\ \Delta_{T_n}W'' 
\right)^2\ind_{A_{n,h}}\\

&=&
\frac{1}{(1+\transp{\varphi}_{T_n}\Delta_{T_n} N)^2}\left(
\transp \left(J''_{T_n}
+
\transp\zeta''_{T_n}{\varphi}_{T_n} \right) \Delta_{T_n}W'' 
\right)^2\\

&=&
\frac{1}{(1+\transp{\varphi}_{T_n}\Delta_{T_n} N)^2}\left(
\Delta_{T_n}D
+
\transp {\varphi}_{T_n}\Delta_{T_n} N 
\right)^2.
\dce
$$
The lemma is proved.
\ok

\textbf{Discussion.}\footnote{A colleague has expressed his belief that Lemma \ref{referee} should have a more direct proof from Lemma \ref{A-G-p}. His comment motivates the present discussion, leading to a quick proof of Lemma \ref{referee}.}
The essential of Lemma \ref{referee} is to compute $\transp K'' \Delta \widetilde{W}''$ for a process $K''$ which satisfies the accessible raw structure condition (\ref{structure-condition-da}). Let $1\leq n<\mathsf{N}^a$. 
Compute the jump at $T_{n}\leq T, T_{n}<\infty,$ of the accessible raw structure condition (\ref{structure-condition-da}). We have$$
\transp K''_{T_{n}} \mathbb{E}[\Delta_{T_{n}} \widetilde{W}''\ \Delta_{T_{n}} \widetilde{W}''_{h} | \mathcal{G}_{T_{n}-}]
=
\mathbb{E}[\Delta_{T_{n}}D\ \transp \Delta_{T_{n}}{W}''_{h} | \mathcal{F}_{T_{n}-}] +\transp{\varphi}_{T_{n}} \mathbb{E}[\Delta_{T_{n}}N\ \transp \Delta_{T_{n}}{W}''_{h} | \mathcal{F}_{T_{n}-}],
$$
for $0\leq h\leq d$. This is equivalent to 
$$
\transp K''_{T_{n}} \frac{1}{2^n}(\epsilon_{h} - \overline{p}_{n}) \overline{p}_{n,h}
=
d_{n,h} p_{n,h} +\transp{\varphi}_{T_{n}} n_{n,h}p_{n,h}.
$$
By Lemma \ref{F-Tn} and Remark \ref{1+fn>0bis}, it is again equivalent to 
$$
\transp K''_{T_{n}} \frac{1}{2^n}(\epsilon_{h} - \overline{p}_{n}) 
=
\frac{d_{n,h}  +\transp{\varphi}_{T_{n}} n_{n,h}}{1 +\transp{\varphi}_{T_{n}} n_{n,h}},
$$
on $\{p_{n,h}>0\}$. But $A_{n,h}\subset \{p_{n,h}>0\}$. We obtain$$
\dcb
\transp K''_{T_{n}} \Delta_{T_{n}} \widetilde{W}'' \ind_{A_{n,h}}&=& \transp K''_{T_{n}} \frac{1}{2^n}(\epsilon_{h} - \overline{p}_{n}) \ind_{A_{n,h}}\\

&=&\frac{d_{n,h}  +\transp{\varphi}_{T_{n}} n_{n,h}}{1 +\transp{\varphi}_{T_{n}} n_{n,h}} \ind_{A_{n,h}}
=
\frac{ \Delta_{T_{n}} D  +\transp{\varphi}_{T_{n}}  \Delta_{T_{n}} N}{1 +\transp{\varphi}_{T_{n}} \Delta_{T_{n}} N } \ind_{A_{n,h}}.
\dce
$$
This proves (for a second time) Lemma \ref{referee}. To end this discussion, we notice that the previous longer proof of Lemma \ref{referee} remains interesting because it deals with filtration changes by probability changes. \ \ok

As a corollary of Lemma \ref{thmprimeprime} and Lemma \ref{referee}, we state

\bethe\label{Rthmprimeprime}
Suppose the drift multiplier assumption \ref{assump1} and Assumption \ref{assump-mrt}. 
The accessible raw structure condition (\ref{structure-condition-da}) related to $D$ is satisfied, if and only if Condition \ref{1+fin} holds and the process $$
\sqrt{\sum_{n=1}^{\mathsf{N}^a-}\ind_{\{T_n\leq t\wedge T\}}\frac{1}{(1+\transp{\varphi}_{T_n}\Delta_{T_n} N)^2}\left(
\Delta_{T_n}D
+
\transp {\varphi}_{T_n}\Delta_{T_n} N  
\right)^2}
$$
is $(\mathbb{P},\mathbb{G})$ locally integrable.
\ethe

\

\section{Solution of the totally inaccessible raw structure condition}

As in the previous section, the drift multiplier assumption \ref{assump1} and Assumption \ref{assump-mrt} are in force.

\subsection{Equations at the stopping times $S_{n}$}

Let $(S_n)_{1\leq n<\mathsf{N}^i}$ ($\mathsf{N}^i\leq \infty$) be a sequence of $(\mathbb{P},\mathbb{F})$ totally inaccessible stopping times such that $[S_n]\cap [S_{n'}]=\emptyset$ for $n\neq n'$ and $\{s\geq 0:\Delta_sW'''\neq 0\}\subset\cup_{n\geq 1}[S_n]$. A $\mathbb{G}$ predictable process $K'''$ satisfies the equation (\ref{structure-condition-di}) if and only if, for $1\leq h<\mathsf{n}'''$, 
\begin{equation}\label{di-GF}
\dcb
K'''_h\centerdot[\widetilde{W}'''_h,\widetilde{W}'''_h]^{\mathbb{G}\cdot p}

=
K'''_h\centerdot[{W}'''_h,{W}'''_h]^{\mathbb{G}\cdot p}
=
(J'''_h+\transp{\varphi}\zeta'''_h)\centerdot [W'''_h, {W}'''_h]^{\mathbb{F}\cdot p}
\dce
\end{equation}
on $[0,T]$.

\bl\label{KequSn}
A $\mathbb{G}$ predictable process $K'''$ satisfies the equation (\ref{structure-condition-di}) if and only if, for $1\leq n<\mathsf{N}^i$, for $1\leq h\leq \mathsf{n}'''$, $K'''_{h,S_n}$ satisfies the equation 
$$
\dcb
&&
(1+\transp{\varphi}_{S_n}\mathtt{R}_n)K'''_{h,S_n}\mathbb{E}[\ind_{\{\Delta_{S_n}W'''_h\neq 0\}}|\mathcal{G}_{S_n-}]
=
(J'''_h+\transp{\varphi}\zeta'''_h)_{S_n}\mathbb{E}[\ind_{\{\Delta_{S_n}W'''_h\neq 0\}}|\mathcal{F}_{S_n-}]
\dce
$$
on $\{S_n\leq T, S_n<\infty\}$, where $\mathtt{R}_n=\mathbb{E}[\Delta_{S_n}N|\mathcal{F}_{S_n-}]$.
\el

\textbf{Proof.}
Let $1\leq n<\mathsf{N}^i, 0\leq h\leq \mathsf{n}'''$. We define $g_{n,h}$ to be an $\mathbb{F}$ (resp. $\overline{g}_{n,h}$ a $\mathbb{G}$) predictable process such that $$
\dcb
g_{n,h,S_n}=\mathbb{E}[(\Delta_{S_n} W'''_h)^2|\mathcal{F}_{S_n-}] \
\mbox{ resp. }
\overline{g}_{n,h,S_n}=\mathbb{E}[(\Delta_{S_n} W'''_h)^2|\mathcal{G}_{S_n-}].
\dce
$$
Let $f$ denote the coefficient of the $(\mathbb{P},\mathbb{F})$ martingale $\Delta_{S_n} W'''_h\ind_{[S_n,\infty)}-(\Delta_{S_n} W'''_h\ind_{[S_n,\infty)})^{\mathbb{F}\cdot p}$ in its martingale representation with respect to $W$. By pathwise orthogonality, $f$ has all components null but one (denoted by $f_{h}$) corresponding to $W'''_h$, and, by Lemma \ref{single-jump}, $f_h$ can be modified to be bounded. We have, on the one hand,$$
\dcb
&&f_{h} K'''_h\centerdot[{W}'''_h,{W}'''_h]^{\mathbb{G}\cdot p}
=
K'''_h\centerdot[f_{h} \centerdot{W}'''_h,{W}'''_h]^{\mathbb{G}\cdot p}
=
K'''_h\centerdot[\transp f \centerdot{W},{W}'''_h]^{\mathbb{G}\cdot p}\\
&=&
K'''_h\centerdot[\Delta_{S_n}{W}'''_h\ind_{[S_n,\infty)}-(\Delta_{S_n}{W}'''_h\ind_{[S_n,\infty)})^{\mathbb{F}\cdot p},{W}'''_h]^{\mathbb{G}\cdot p}\\

&=&
K'''_h\centerdot((\Delta_{S_n}{W}'''_h)^2\ind_{[S_h,\infty)})^{\mathbb{G}\cdot p}
=
K'''_h\overline{g}_{n,h}\centerdot\left(\ind_{[S_h,\infty)}\right)^{\mathbb{G}\cdot p}
\dce
$$
on $[0,T]$. On the other hand, 
$$
\dcb
&&
f_{h}(J'''_h+\transp{\varphi}\zeta'''_h)\centerdot [W'''_h, {W}'''_h]^{\mathbb{F}\cdot p}
=
(J'''_h+\transp{\varphi}\zeta'''_h)g_{n,h}\centerdot\left(\ind_{[S_n,\infty)}\right)^{\mathbb{F}\cdot p}
\dce
$$
on $[0,T]$. All put together, equation (\ref{di-GF}) implies
\begin{equation}\label{di-GF2}
K'''_h\overline{g}_{n,h}\centerdot\left(\ind_{[S_n,\infty)}\right)^{\mathbb{G}\cdot p}
=
(J'''_h+\transp{\varphi}\zeta'''_h)g_{n,h}\centerdot\left(\ind_{[S_n,\infty)}\right)^{\mathbb{F}\cdot p}
\end{equation}
on $[0,T]$ for any $1\leq n<\mathsf{N}^i, 0\leq h\leq \mathsf{n}'''$. Conversely, we note that$$
\dcb
&&[{W}'''_h,{W}'''_h]^{\mathbb{F}\cdot p}
=
\sum_{n=1}^{\mathsf{N}^a-}\left((\Delta_{S_n}{W}'''_h)^2\ind_{[S_n,\infty)}\right)^{\mathbb{F}\cdot p}
=
\sum_{n=1}^{\mathsf{N}^a-}g_{n,h}\centerdot\left(\ind_{[S_n,\infty)}\right)^{\mathbb{F}\cdot p},
\dce
$$
and
$$
\dcb
&&[{W}'''_h,{W}'''_h]^{\mathbb{G}\cdot p}
=
\sum_{n=1}^{\mathsf{N}^a-}\left((\Delta_{S_n}{W}'''_h)^2\ind_{[S_n,\infty)}\right)^{\mathbb{G}\cdot p}
=
\sum_{n=1}^{\mathsf{N}^a-}\overline{g}_{n,h}\centerdot\left(\ind_{[S_n,\infty)}\right)^{\mathbb{G}\cdot p}.
\dce
$$
Consequently, if the process $K'''$ satisfied all equation (\ref{di-GF2}) for $1\leq n<\mathsf{N}^i, 0\leq h\leq \mathsf{n}'''$, the process $K'''$ is $[{W}'''_h,{W}'''_h]^{\mathbb{G}\cdot p}$-integrable and equation (\ref{di-GF}) is satisfied.

Consider the equations (\ref{di-GF2}). Following Lemma \ref{A-G-p}, we compute on $[0,T]$ :
\begin{equation}\label{1+fNGFS}
\dcb
&&
(\ind_{[S_n,\infty)})^{\mathbb{G}\cdot p}
=
(\ind_{[S_n,\infty)})^{\mathbb{F}\cdot p}+\transp{\varphi}_{S_n}(\Delta_{S_n}N\ind_{[S_n,\infty)})^{\mathbb{F}\cdot p}
=
\left(1+\transp{\varphi}\mathtt{r}_n\right)\centerdot(\ind_{[S_n,\infty)})^{\mathbb{F}\cdot p},
\dce
\end{equation}
where $\mathtt{r}_n$ is an $\mathbb{F}$ predictable process such that $(\mathsf{r}_n)_{S_n}=\mathtt{R}_n$. Hence, on $[0,T]$, for any $\mathbb{G}$ predictable set $\mathtt{A}$ such that $\ind_{\mathtt{A}}(1+\transp{\varphi}\mathtt{r}_n)$ and $J'''_h+\transp{\varphi}\zeta'''_h$ are  bounded, equation (\ref{di-GF2}) implies
$$
\dcb
&&\ind_{\mathtt{A}}(1+\transp{\varphi}\mathtt{r}_n)
K'''_h\overline{g}_{n,h}\centerdot\left(\ind_{[S_k,\infty)}\right)^{\mathbb{G}\cdot p}
=
\ind_{\mathtt{A}}(1+\transp{\varphi}\mathtt{r}_n)(J'''_h+\transp{\varphi}\zeta'''_h)g_{n,h}\centerdot (\ind_{[S_n,\infty)})^{\mathbb{F}\cdot p}\\

&=&
\ind_{\mathtt{A}}(J'''_h+\transp{\varphi}\zeta'''_h)g_{n,h}\centerdot (\ind_{[S_n,\infty)})^{\mathbb{G}\cdot p}.
\dce
$$
This is equivalent to
$$
\dcb
&&(1+\transp{\varphi}\mathtt{R}_n)K'''_{h,S_n}\overline{g}_{h,n,S_n}
=
(J'''_h+\transp{\varphi}\zeta'''_h)_{S_n}g_{n,h,S_n},\ \mbox{ on $\{S_n\leq T, S_n<\infty\}$}.
\dce
$$
Let $\alpha'''$ be the $\mathbb{F}$ predictable process in subsection \ref{www} such that, for $1\leq h\leq \mathsf{n}'''$, $$
\Delta W'''_h = \alpha'''_h\ind_{\{\Delta W'''_h\neq 0\}}.
$$ 
We compute
$$
\dcb
&&(\alpha'''_{h,S_n})^2
(1+\transp{\varphi}_{S_n}\mathtt{R}_n)K'''_{h,S_n}\mathbb{E}[\ind_{\{\Delta_{S_n}W'''_h\neq 0\}}|\mathcal{G}_{S_n-}]
=
(1+\transp{\varphi}_{S_n}\mathtt{R}_n)K'''_{h,S_n}\mathbb{E}[(\Delta_{S_n}W'''_h)^2|\mathcal{G}_{S_n-}]\\
&=&
( J'''_h+\transp{\varphi}\zeta'''_h)_{S_n}\mathbb{E}[(\Delta_{S_n}W'''_h)^2|\mathcal{F}_{S_n-}]
=
(\alpha'''_{h,S_n})^2( J'''_h+\transp{\varphi}\zeta'''_h)_{S_n}\mathbb{E}[\ind_{\{\Delta_{S_n}W'''_h\neq 0\}}|\mathcal{F}_{S_n-}].  
\dce
$$
The lemma is proved, because on $\{\alpha'''_{h,S_n}=0\}$, $$
\mathbb{E}[\ind_{\{\Delta_{S_n}W'''_h\neq 0\}}|\mathcal{G}_{S_n-}]
=
\mathbb{E}[\ind_{\{\Delta_{S_n}W'''_h\neq 0\}}|\mathcal{F}_{S_n-}]
=0.\ \ok
$$

\

\subsection{Conditional expectations at stopping times $S_n$}\label{exp-at-jump}

For a fixed $1\leq n<\mathsf{N}^i$, applying the martingale representation property, applying \cite[Lemme(4.48)]{Jacodlivre} with the finite $\mathbb{F}$ predictable constraint condition in subsection \ref{www}, we see that, on $\{S_n<\infty\}$,$$
\mathcal{F}_{S_n}=\mathcal{F}_{S_n-}\vee\sigma(\Delta_{S_n}W''')
=\mathcal{F}_{S_n-}\vee\sigma(\{\Delta_{S_n}W'''_1\neq 0\},\ldots,\{\Delta_{S_n}W'''_{\mathsf{n}'''}\neq 0\}).
$$
(Note that $\{S_n<\infty\}\subset\{\Delta_{S_n}W'''\neq 0\}$.) We set$$
B_{n,k}:=\{\Delta_{S_n}W'''_k\neq 0\},\ q_{n,k}=\mathbb{P}[B_{n,k}|\mathcal{F}_{S_n-}],\ 
\overline{q}_{n,k}=\mathbb{P}[B_{n,k}|\mathcal{G}_{S_n-}],\ 1\leq k\leq {\mathsf{n}'''}.
$$
Note that, by our choice of $W'''$ in subsection \ref{www}, the $B_{n,k}$ form a partition on $\{S_n<\infty\}$ (cf. \cite{song-mrp-drift}).

\bl\label{F-Sn}
Let $1\leq n<\mathsf{N}^i$.   
\ebe
\item[. ]
For any finite random variable $\xi\in\mathcal{F}_{S_n}$, the conditional expectation $\mathbb{E}[\xi|\mathcal{F}_{S_n-}]$ is well-defined. Let $$
\mathfrak{i}_n(\xi)_k=\ind_{\{q_{n,k}>0\}}\frac{1}{q_{n,k}}\mathbb{E}[\ind_{B_{n,k}}\xi|\mathcal{F}_{S_n-}], \ 1\leq k\leq \mathsf{n}'''.
$$
We have $\xi=\sum_{h=1}^{\mathsf{n}'''}\mathfrak{i}_n(\xi)_h\ind_{B_h}$. 
 
\item[. ]
Denote the vector valued random variable $l_{n,k}:=\mathfrak{i}_n(\Delta_{S_n}N)_k, 1\leq k\leq {\mathsf{n}'''}$. We have $$
(1+\transp{\varphi}_{S_n}l_{n,k})q_{n,k}
=
(1+\transp{\varphi}_{S_n}\mathtt{R}_{n})\mathbb{E}[\ind_{B_{n,k}}|\mathcal{G}_{S_n-}]
=
(1+\transp{\varphi}_{S_n}\mathtt{R}_{n})\overline{q}_{n,k}
$$
for $1\leq k\leq {\mathsf{n}'''}$, on $\{S_n\leq T, S_n<\infty\}$, where $\mathtt{R}_n$ is the vector valued process introduced in Lemma \ref{KequSn}.

\item[. ]
We have $(1+\transp{\varphi}_{S_n}\mathtt{R}_{n})>0$ almost surely on $\{S_n\leq T, S_n<\infty\}$.

\dbe

\el

\textbf{Proof.}
The proof of the first assertion is straightforward. To prove the second assertion, we introduce $\mathbb{F}$ predictable processes $H,G$ such that $H_{S_n}=l_{n,k}, G_{S_n}=q_{n,k}$. We apply then Lemma \ref{A-G-p} to write, for any $1\leq k\leq d$, $$
\dcb
(\ind_{B_{n,k}}\ind_{[S_n,\infty)})^{\mathbb{G}\cdot p}

&=&(1+\transp{\varphi}H)G\centerdot(\ind_{[S_n,\infty)}
)^{\mathbb{F}\cdot p}
\dce
$$
on $[0,T]$. Integrate the term $(1+\transp{\varphi}\mathtt{r}_{n})$ and apply formula (\ref{1+fNGFS}) to obtain
$$
\dcb
((1+\transp{\varphi}_{S_n}\mathtt{R}_{n})\ind_{B_k}\ind_{[S_n,\infty)})^{\mathbb{G}\cdot p}

&=&((1+\transp{\varphi}_{S_n}l_{n,k})q_{n,k}\ind_{[S_n,\infty)}
)^{\mathbb{G}\cdot p}\\

\dce
$$
on $[0,T]$, which proves the second formula. Consider $(1+\transp{\varphi}_{S_n}\mathtt{R}_{n})$. We compute on $[0,T]$ $$
\dcb
0\leq (\ind_{\{1+\transp{\varphi}_{S_n}\mathtt{R}_{n}\leq 0\}}\ind_{[S_n,\infty)})^{\mathbb{G}\cdot p}

&=&\ind_{\{1+\transp{\varphi}\mathtt{r}_{n}\leq 0\}}(1+\transp{\varphi}\mathtt{r}_{n})\centerdot(\ind_{[S_n,\infty)})^{\mathbb{F}-p}\leq 0.
\dce
$$
This yields $
\mathbb{E}[\ind_{\{1+\transp{\varphi}_{S_n}\mathtt{R}_{n}\leq 0\}}\ind_{\{S_n\leq T,S_n<\infty\}}]
=0,
$
proving the third assertion.  \ok

\bl\label{1+fn2>0}
For $ 1\leq n<\mathsf{N}^i$ we have $
1+\transp{\varphi}_{S_n}\Delta_{S_n} N> 0
$
on $\{S_n\leq T, S_n<\infty\}$.
\el

\textbf{Proof.} We compute, for $1\leq h\leq \mathsf{n}'''$,$$
\dcb
&&0\leq\mathbb{E}[\ind_{\{1+\transp{\varphi}_{S_n}\Delta_{S_n} N \leq 0\}}\ind_{B_{n,h}}|\mathcal{G}_{S_n-}]

=\ind_{\{1+\transp{\varphi}_{S_n}l_{n,h} \leq 0\}}\overline{q}_{n,h}
\\

&=&\ind_{\{1+\transp{\varphi}_{S_n}l_{n,h} \leq 0\}}\frac{1+\transp{\varphi}_{S_n}l_{n,h}}{1+\transp{\varphi}_{S_n}\mathtt{R}_n}q_{n,h}\leq 0
\dce
$$
 on $\{S_n\leq T, S_n<\infty\}$. It follows that $\ind_{\{1+\transp{\varphi}_{S_n}\Delta_{S_n} N \leq 0\}}\ind_{B_{n,h}}=0$ for $1\leq h\leq \mathsf{n}'''$, i.e., $1+\transp{\varphi}_{S_n}\Delta_{S_n} N>0$. \ok

\

\subsection{Consequences on the totally inaccessible raw structure condition}

Note that $\widetilde{W}'''$ is a $\mathbb{G}$ purely discontinuous local martingale. This is because $W'''$ is the limit in martingale space of $(\mathbb{P},\mathbb{F})$ local martingales with finite variation (cf. \cite[Theorem 6.22]{HWY}), and therefore the same is true for $\widetilde{W}'''$ by \cite[Proposition (2,2)]{J}.

\bethe\label{KWG2}
Suppose the drift multiplier assumption \ref{assump1} and Assumption \ref{assump-mrt}. 
The totally inaccessible raw structure condition (\ref{structure-condition-di}) related to $D$ is satisfied for all $1\leq h\leq \mathsf{n}'''$, if and only if the process
$$
\sqrt{\sum_{n=1}^{\mathsf{N}^i-}\ind_{\{S_n\leq t\wedge T\}}\frac{1}{(1+\transp{\varphi}_{S_n}\Delta_{S_n} N)^2}\left(
\Delta_{S_n}D
+
\transp {\varphi}_{S_n}\Delta_{S_n} N  
\right)^2}
$$
is $(\mathbb{P},\mathbb{G})$ locally integrable. In this case, a solution process $K'''$ is given by
\begin{equation}\label{sol-i}
\dcb
&&
K'''_{h}
=
\frac{J'''_h+\transp{\varphi}\zeta'''_h}{1+\transp{\varphi}\zeta'''_h\alpha'''_h}\ind_{\{1+\transp{\varphi}\zeta'''_h\alpha'''_h\neq 0\}},\ \mbox{ $\mathsf{d}[W'''_h,W'''_h]-a.s.$ on $[0,T]$, $1\leq h\leq \mathsf{n}'''$.} 
\dce
\end{equation}
\ethe

\textbf{Proof.} 
Suppose the integrability condition in the theorem and define $K'''_h$, $1\leq h\leq \mathsf{n}'''$, by (\ref{sol-i}). As$$
l_{n,h}\ind_{B_{n,h}}
=
\Delta_{S_n}N\ind_{B_{n,h}}
=
\zeta_{S_n}\Delta_{S_n}W\ind_{B_{n,h}}
=
\zeta'''_{h,S_n}\Delta_{S_n}W'''_h\ind_{B_{n,h}}
=
(\zeta'''_h\alpha'''_h)_{S_n}\ind_{B_{n,h}},
$$
the formula (\ref{sol-i}) implies, on $\{S_n\leq T, S_n<\infty\}$ for any $1\leq n<\mathsf{N}^i$, $$
\dcb
&&
K'''_{h,S_n}\ind_{B_{n,h}}
=
\frac{(J'''_h+\transp{\varphi}\zeta'''_h)_{S_n}}{(1+\transp{\varphi}_{S_n}l_{n,h})}\ind_{B_{n,h}}.
\dce
$$
(Noting that the random measure $\mathsf{d}[W'''_h,W'''_h]$ charges the set $B_{n,h}\cap[S_n]$). Take the conditioning with respect to $\mathcal{G}_{S_n-}$ with help of Lemma \ref{F-Sn}.
\begin{equation}\label{lrqK}
\dcb
&&
(1+\transp{\varphi}_{S_n}l_{n,k})q_{n,k} K'''_{h,S_n}
=
(J'''_h+\transp{\varphi}\zeta'''_h)_{S_n}q_{n,k}, 
\dce
\end{equation}
i.e., the equations in Lemma \ref{KequSn} are satisfied. We prove hence that $K'''_h$ is a solution of equation (\ref{structure-condition-di}). 

We now prove that $K'''_h$ is $\widetilde{W}'''_h$-integrable on $[0,T]$. For any $1\leq n<\mathsf{N}^i$, on the set $B_{n,h}\cap\{S_n\leq T, S_n<\infty\}$, $J'''_h\Delta_{S_n}W'''_h=\Delta_{S_n}D$, $(\transp{\varphi}\zeta'''_h)_{S_n}\Delta_{S_n}W'''_h=\transp{\varphi}_{S_n}\Delta_{S_n} N$ so that 
\begin{equation}\label{KjumpW}
\dcb
K'''_{h,S_n}\Delta_{S_n}W'''_h
&=&\frac{(J'''_h+\transp{\varphi}\zeta'''_h)_{S_n}}{(1+\transp{\varphi}\zeta'''_h\alpha'''_h)_{S_n}}\Delta_{S_n}W'''_h\\
&=&
\frac{J'''_h\Delta_{S_n}W'''_h+(\transp{\varphi}\zeta'''_{h})_{S_n}\Delta_{S_n}W'''_h}{1+(\transp{\varphi}\zeta'''_h)_{S_n}\Delta_{S_n}W'''_h}
=
\frac{\Delta_{S_n}D+\transp{\varphi}\Delta_{S_n} N}{1+\transp{\varphi}\Delta_{S_n} N}\ind_{B_{n,h}}.
\dce
\end{equation}
This proves the $\widetilde{W}'''_h$-integrability of $K'''_h$ on $[0,T]$. (Recall $\Delta\widetilde{W}'''=\Delta{W}'''$.) 

We finally check if $K'''_h\Delta \widetilde{W}'''_h=K'''_h\Delta W'''_h<1$ on $[0,T]$. But, on $B_{n,h}\cap\{S_n\leq T, S_n<\infty\}$,
$$
\dcb
&&
K'''_{h,S_n}\Delta_{S_n}W'''_h
=
\frac{\Delta_{S_n}D+\transp{\varphi}\Delta_{S_n} N}{1+\transp{\varphi}\Delta_{S_n} N}
\

<
\frac{1+\transp{\varphi}\Delta_{S_n} N}{1+\transp{\varphi}\Delta_{S_n} N} =1,
\dce
$$
because $1+\transp{\varphi}\Delta_{S_n} N>0$ (cf. Lemma \ref{1+fn2>0}). The totally inaccessible raw structure condition related to $D$ is satisfied by $K'''_h$, $1\leq h\leq \mathsf{n}'''$.

Conversely, suppose that $K'''_h$ is a solution of the totally inaccessible raw structure condition (\ref{structure-condition-di}) related to $D$. The formula in Lemma \ref{KequSn} is satisfied so as formula (\ref{lrqK}) (with help of Lemma \ref{F-Sn}). Multiply formula (\ref{lrqK}) by $\ind_{B_{n,h}}$ on $\{S_n\leq T, S_n<\infty\}$, we obtain
$$
\dcb
&&
K'''_{h,S_n}\ind_{B_{n,k}}
=
\frac{(J'''_h+\transp{\varphi}\zeta'''_h)_{S_n}}{1+\transp{\varphi}_{S_n}l_{n,h}}\ind_{B_{n,k}}
=
\frac{(J'''_h+\transp{\varphi}\zeta'''_h)_{S_n}}{(1+\transp{\varphi}\zeta'''_h\alpha'''_h)_{S_n}}\ind_{B_{n,k}},
\dce
$$
for $1\leq n<\mathsf{N}^i$, $1\leq h\leq \mathsf{n}'''$. This implies in turn the formula (\ref{KjumpW}). The $\widetilde{W}'''_h$-integrability of $K'''_h$ on $[0,T]$ for $1\leq h\leq \mathsf{n}'''$ implies finally the integrability condition of the theorem. \ok

\

\section{Proof of Theorems 4.2 and 4.3}

\textbf{Proof of Theorem \ref{main}}
It is direct consequence of Lemma \ref{piece-ensemble} together with Theorem \ref{answer-c}, Theorem \ref{Rthmprimeprime} and Theorem \ref{KWG2}. \ok

\

\textbf{Proof of Theorem \ref{fullviability}}
We have the martingale representation property with representation process $W=(W',W'',W''')$ in Assumption \ref{assump-mrt}. 

\textbf{Necessary part}

Suppose the full viability on $[0,T]$. Then, the drift operator satisfied the drift multiplier assumption on $[0,T]$, as it is proved in \cite[Theorem 5.5]{song-mrp-drift}. Let us prove that the group of the conditions (\ref{structure-condition-c}), (\ref{structure-condition-da}), (\ref{structure-condition-di}) related to $D=0$ in $\mathbb{G}$ are satisfied.

First of all, there exists a $\mathbb{G}$ structure connector $Y$ for $W'$ on $[0,T]$ (cf. Remark \ref{locallyboundedM}), i.e. (cf. formula (\ref{structure-condition}))
$$
[Y,\widetilde{W}']^{\mathbb{G}\cdot p}=\transp{\varphi}\centerdot [N,W']^{\mathbb{F}\cdot p}.
$$
By the continuity, we can replace $[Y,\transp\widetilde{W}']^{\mathbb{G}\cdot p}$ by $\transp K'\centerdot[\widetilde{W}',\transp \widetilde{W}']^{\mathbb{G}\cdot p}$ for some $\mathbb{G}$ predictable $\widetilde{W}'$-integrable process $K'$. We prove thus the continuous raw structure condition (\ref{structure-condition-c}) related to $D=0$ in $\mathbb{G}$.

Let us check that the conditions in \cite[Theorem 3.8]{song-mrp-drift} is satisfied for the jump measure of $\widetilde{W}''$. Recall $\boldsymbol{\epsilon}=(\epsilon_h)_{0\leq h\leq d}$ the canonical basis in $\mathbb{R}\times\mathbb{R}^{d}$, $\mathsf{a}_{n}$ the vector $(\ind_{A_{n,h}})_{0\leq h\leq d}$, $1\leq n<\mathsf{N}^a$. We have the identity on $\{T_n<\infty\}$
$$
\Delta_{T_n}\widetilde{W}''=\frac{1}{2^n}(\mathsf{a}_{n}-\overline{p}_n)
=
\frac{1}{2^n}\sum_{h=0}^d(\epsilon_h-\overline{p}_n)\ind_{A_{n,h}}.
$$
Conforming to the notation of \cite[Theorem 3.8]{song-mrp-drift}, we have $$
\alpha_{n,h}=\frac{1}{2^n}(\epsilon_h-\overline{p}_n), \ \
\gamma_{n,i}= \frac{1}{2^n}(\epsilon_i-\overline{p}_{n,i}\mathbf{1}), \
0\leq i\leq d.
$$
where $\transp \mathbf{1}=(1,1,\ldots,1)\in \mathbb{R}\times\mathbb{R}^{d}$. Let $v$ be a vector in $\mathbb{R}\times\mathbb{R}^{d}$ orthogonal to the $\gamma_{n,i}$. Then$$
v_i=\overline{p}_{n,i}\transp \mathbf{1}v,\ 0\leq i\leq d,
$$
i.e., $v$ is proportional to $\overline{p}_n$. The vectors $\gamma_i$ together with $\overline{p}_n$ span whole $\mathbb{R}\times\mathbb{R}^{d}$, which is the condition of \cite[Theorem 3.8]{song-mrp-drift}.

By the full viability, there exists a $\mathbb{G}$ structure connector $Y$ for $W''$ on $[0,T]$, i.e. 
$$
[Y,\widetilde{W}'']^{\mathbb{G}\cdot p}=\transp{\varphi}\centerdot [N,W'']^{\mathbb{F}\cdot p}.
$$
Applying \cite[Lemma 3.1 and Theorem 3.8]{song-mrp-drift}, we can replace $[Y,\transp\widetilde{W}'']^{\mathbb{G}\cdot p}$ by $\transp K''\centerdot[\widetilde{W}'',\transp \widetilde{W}'']^{\mathbb{G}\cdot p}$ for some $\mathbb{G}$ predictable $\widetilde{W}''$-integrable process $K''$. Moreover, the time support of the jump measure of $\widetilde{W}''$ is $\cup_{1\leq n<\mathsf{N}^a}[T_{n}]$ so that the process $\widehat{U}$ in \cite[Lemma 3.1]{song-mrp-drift} is null, according to \cite[Theorem (3.75) (a)]{Jacodlivre}. Following \cite[Lemma 3.1 and Theorem 3.8]{song-mrp-drift}, we see that $\transp K''\Delta\widetilde{W}''$ is the conditional expectation of $\Delta Y$ given the $\sigma$-algebra $\widetilde{\mathcal{P}}$ under the Dolean-Dade measure associated with the jump of $\widetilde{W}''$. Hence, $\transp K''\Delta\widetilde{W}''<1$, proving the accessible raw structure condition (\ref{structure-condition-da}) related to $D=0$ in $\mathbb{G}$.

Notice that $\Delta \widetilde{W}'''_h$, $1\leq h\leq \mathsf{n}'''$, satisfies clearly the condition in \cite[Theorem 3.9]{song-mrp-drift} (with $\mathsf{n}=1$). We can repeat the above reasoning to prove the totally inaccessible raw structure condition (\ref{structure-condition-di}) related to $D=0$ in $\mathbb{G}$, noting that the time support of the jump measure of $\widetilde{W}'''$ is $\cup_{1\leq n<\mathsf{N}^i}[S_{n}]$ so that the process $\widehat{U}$ in \cite[Lemma 3.1]{song-mrp-drift} is null.

Now apply Theorem \ref{answer-c}, Theorem \ref{Rthmprimeprime} and Theorem \ref{KWG2}. We prove Condition \ref{1+fin} and the condition (\ref{fn-sur-fn}) in Theorem \ref{fullviability}.

\textbf{Sufficient part}

Conversely, suppose the drift multiplier assumption, Condition \ref{1+fin} and the condition (\ref{fn-sur-fn}). Apply Lemma \ref{piece-ensemble} together with Theorem \ref{answer-c}, Theorem \ref{Rthmprimeprime} and Theorem \ref{KWG2}, then translate the conclusion with Theorem \ref{deflator-connector} in term of deflators. We conclude that any $\mathbb{F}$ local martingale has a $\mathbb{G}$ deflator. Let now $S$ be an $\mathbb{F}$ special semimartingale with a $\mathbb{F}$ deflator $D$. Then, $(D,DS)$ has a $\mathbb{G}$ deflator $Y$, i.e. $S$ has a $\mathbb{G}$ deflator $DY$. Apply again Theorem \ref{deflator-connector}. We conclude the proof.\ \ok

\

\textbf{Proof of Corollary \ref{commondeflator}}
We note that, in the case $D\equiv 0$, the proof of Lemma \ref{piece-ensemble} gives a common $(\mathbb{P},\mathbb{G})$ structure connector to all $(\mathbb{P},\mathbb{F})$ local martingales. Corollary \ref{commondeflator} is therefore the consequence of Theorem \ref{deflator-connector}. \ok

\

\section{Theorem \ref{fullviability} in applications}

We underline that the most important aspect of Theorem \ref{fullviability}  is its theoretical conclusion about the drift multiplier assumption and Condition \ref{1+fin} in relation with the full viability problem. Theorem \ref{fullviability} is considered, because we need to know what exactly are the factors which determine the full viability in an expanded information flow, and we need to know what we can say and what we can do, if the full viability is expected.

The drift multiplier assumption and Condition \ref{1+fin} introduced in section \ref{DMA} are not very common in the literature of enlargement of filtration. Condition \ref{1+fin} states that the conditional probability measures on the $\sigma$-algebra $\mathcal{F}_{R}$, given $\mathcal{F}_{R-}$ or given $\mathcal{G}_{R-}$ are equivalent. We have already indicated in Remark \ref{rq:cR} that the quasi-left-continuity of $\mathbb{F}$ implies Condition \ref{1+fin}. We can also understand this condition in the following way. Suppose that there exist two random variables $\xi$ and $\zeta$ such that $$
\mathcal{F}_{R} = \mathcal{F}_{R-}\vee\sigma(\xi), \ \ \mathcal{G}_{R-} = \mathcal{F}_{R-}\vee\sigma(\zeta).
$$
Consider every thing under the conditional probability $\mathbb{P}^{:}=\mathbb{P}[\cdot\ |\mathcal{F}_{R-}]$ so that the $\sigma$-algebra $\mathcal{F}_{R-}$ becomes $\mathbb{P}^{:}$-trivial. Then, Condition \ref{1+fin} means that the conditional law of $\xi$ given $\sigma(\zeta)$ under $\mathbb{P}^{:}$ is equivalent to its unconditional law under $\mathbb{P}^{:}$, or, in other words, $\xi$ and $\zeta$ can become independent by an equivalent change of the probability measure $\mathbb{P}^{:}$.

The drift multiplier assumption in Definition \ref{assump1} is a necessary condition, whenever the martingale representation property for $\mathbb{F}$ and the full viability for $\mathbb{G}$ are satisfied, as explained in \cite[Theorem 5.5]{song-mrp-drift}. However, the idea of drift multiplier assumption is inspired from the classical models: the initial enlargement of filtration with a random variable $\xi$ (valued in a measurable space $(\mathtt{E},\mathcal{B})$) under Jacod criterion in \cite{Jacod}, or the progressive enlargement with a general random time $\tau$ restricted on the time horizon $[0,\tau]$ in \cite{J, JY2}, or the progressive enlargement with a honest time $\tau$ in \cite{barlow, J, JY2}, or the progressive enlargement with an initial time in \cite{CJZ, JC}. According to \cite[Theorem (2.1)]{Jacod} (with the notations therein), under Jacod criterion, the drift operator is given by $$
\Gamma(X) = \left(\frac{1}{q^x_{-}}\centerdot\cro{q^x,X}^{\mathbb{P}\cdot\mathbb{F}}\right)_{x=\xi},
$$
for any $\mathbb{F}$ local martingale $X$, where $q^x$ denote the density functions introduced in \cite{Jacod}. If the martingale representation property holds in $\mathbb{F}$ with a locally bounded representation process $W$, the oblique bracket takes the form$$
\cro{q^x,X}^{\mathbb{P}\cdot\mathbb{F}} = \transp H^x\centerdot\cro{W,X}^{\mathbb{P}\cdot\mathbb{F}},
$$
for some $\mathcal{B}\otimes\mathcal{P}(\mathbb{F})$ measurable vector valued function $H^x_{t}(\omega)$. Using the notations in Definition \ref{assump1}, this means that the drift multiplier assumption is satisfied with $$
\varphi = \frac{H^\xi}{q^\xi_{-}}, \ \ N = W, \ \ \transp\varphi\Delta N = \frac{\transp H^\xi \Delta W}{q^\xi_{-}}
= \frac{\Delta q^\xi}{q^\xi_{-}} \ \mbox{ and }\ \frac{\transp\varphi\Delta N}{1+\transp\varphi\Delta N} = \frac{\Delta q^\xi}{q^\xi_{}}.
$$ 
The condition (\ref{fn-sur-fn}) in Theorem \ref{fullviability} has now a natural interpretation in term of the density functions $q^x$ of \cite{Jacod} (more precisely, in term of $\frac{1}{(q^\xi)^2}\centerdot [q^x,q^x]|_{x=\xi}$).

Consider a progressive enlargement with a random time $\tau$. To make the presentation easier, assume that $0<\tau<\infty$ and $\tau$ avoids all $\mathbb{F}$ stopping times. Then,  on the horizon $[0,\tau]$, according to the drift formula in \cite{J, JY2}, the process $N$ in the drift multiplier assumption can be the martingale part of the Az\'ema supermartingale $Z$ (in the filtration $\mathbb{F}$) of $\tau$ and $\varphi=\frac{1}{Z_{-}}$. We see then $\frac{\transp\varphi\Delta N}{1+\transp\varphi\Delta N} = \frac{\Delta Z}{Z}$ on $[0,\tau]$, which gives a good interpretation of the condition (\ref{fn-sur-fn}) in term of $Z$ (more precisely, in term of $\frac{1}{Z^2}\centerdot [Z,Z]$).

For the applications, being necessary and sufficient condition, Theorem \ref{fullviability} can be applied in the two senses. We present below an example where Theorem \ref{fullviability} is used to construct fully viable multi-default time models. 

We have the following idea on the incomplete markets. Initially, every market is fair, smooth and complete like the Black-Scholes model. Only over the time, successive default events and crisis deteriorate the market condition. The market persists, but becomes more and more unpredictable and incomplete. A natural way to model such a situation is to begin with a fair filtration $\mathbb{F}$, and then to expand $\mathbb{F}$ successively with random times $\tau_{1},\ldots, \tau_{n}$ (multi-default time model), where the random times $\tau_{1},\ldots, \tau_{n}$ must be chosen to preserve the viability of the market.  

The most important point, in the construction of fully viable multi-default time models with Theorem \ref{fullviability}, is that, at every step of the successive enlargements, the three properties should be satisfied, namely the martingale representation property, the drift multiplier assumption and Condition \ref{1+fin} (or the quasi-left-continuity). To establish the three properties, one should know the enlargement of filtration formulas (on the whole time horizon $\mathbb{R}_{+}$) at every step of the successive enlargements. That is a problem, because no general formula is known for the successive enlargements with general random times $\tau_{1},\ldots, \tau_{n}$. Honest times constitute a class of random times for which enlargement of filtration formulas exist (cf.\cite{J}). However, it is not a good choice for market modeling, because the honest times typically create arbitrage opportunity (cf. \cite{FJS}). The class of initial times (i.e., times satisfying Jacod's criterion) introduced in \cite{JC} are good candidates to do successive enlargement as in \cite{pham}. The three properties can be established if the density functions $q^x$ are good enough. There is a disadvantage to work with initial time models, because it is difficult to calibrate the density functions $q^x$ with market data. A third class of random times, for which enlargement of filtration formulas exist, is the class of random times issued from $\natural$-models, introduced in \cite{JS2, songmodel}. The advantage of $\natural$-models is that one can define directly the drift operators. It is a useful property when the models are to be calibrated with market data. A complete study of fully viable multi-default time $\natural$-models is made in the paper \cite{songsuccessive} with Theorem \ref{fullviability}. 

To end this section, we mention Corollary \ref{commondeflator} which implies that, under the martingale representation property, for any fully viable enlargement of filtration model, the enlargement of filtration formula will be necessarily a generalized Girsanov formula, according to the local solution method developed in \cite{song, song-l-s}.

\

\

\appendix

\section{An error in the previous version}

With respect to the previous version, Condition \ref{1+fin} is added in this new version. Here is the reason.

The basic question is whether the condition (\ref{ptop}) of Lemma \ref{samekernels} holds.
 In the previous version, this condition was taken for granted, because of the following passage (at the end of the old proof).
\begin{quote}
$(**): $\texttt{"}
By Lemma 7.1 and Lemma 7.2, $\{0\leq h\leq d: p_{n,h}>0\}=\{0\leq h\leq d: \overline{p}_{n,h}>0\}$.\texttt{"}
\end{quote}
The reasoning $(**)$ was based on the following implicit argument. \begin{quote}
\texttt{"}
The positivity $
1+\transp{\varphi}_{T_n}\Delta_{T_n} N> 0
$ in Lemma 7.2 implies the positivity $(1+\transp{\varphi}_{T_n}n_{n,h})>0$.\texttt{"}
\end{quote}
However, this argumentation leaves a gap. We only know that the random variable $n_{n,h}$ coincides with $\Delta_{T_{n}}N$ on the set $A_{n,h}$. We have no information about $n_{n,h}$ outside of $A_{n,h}$. Hence, the above implication can not be checked outside of $A_{n,h}$.

This is why, in this version, the condition (\ref{ptop}) (under the form of Condition \ref{1+fin}) is inserted as an element in the balance of the necessary and sufficient conditions for the accessible raw structure condition (\ref{structure-condition-da}) related to $D$ to have a solution.

\end{document}